\documentclass[reqno,ammssymb]{article}
\usepackage{bbm}
\usepackage{color}
\usepackage{xcolor}
\usepackage{mathtools}
\usepackage{mathrsfs}
\usepackage{amsfonts}
\usepackage{dsfont}
\usepackage{amscd}
\usepackage{diagbox}
\usepackage{amsthm,amsmath}
\usepackage{xypic}
\usepackage{rotating}
\usepackage{arydshln}
\usepackage[all]{xy}
\usepackage{amssymb}
\usepackage{extarrows}
\usepackage{multirow}
\usepackage{mathdots}
\usepackage{xypic}
\usepackage{ragged2e}
\usepackage{cite}
\usepackage{kpfonts} 

\numberwithin{equation}{section}
\usepackage[top=25.4mm, bottom=25.4mm, left=36mm, right=36mm]{geometry}

\def\otm{\otimes}

\def\Z{\mathbb Z}
\def\K{\mathbb{K}}

\def\lr#1{\left\langle\,#1\right\rangle}

\def\Z{{\mathbb Z}}

\def\ad{{\cdot_{\rm ad}}}
\def\cs#1{{\rm (see \cite{#1})}}
 
\def\noa{\begin{notation}}
\def\hyp{\begin{hypothesis}}
\def\lbb#1{\left[\kern-0.1em\left[ {#1}
 \right]\kern-0.1em\right]}

\newtheorem{thm}{Theorem}[section]
\newtheorem{lem}[thm]{Lemma}
\newtheorem{prop}[thm]{Proposition}

\newtheorem{defn}[thm]{Definition}

\newtheorem{exam}[thm]{Example}

\allowdisplaybreaks[4]
\newcommand{\affil}[2][1]{\textsuperscript{#1}#2}
\date{}
\begin{document}

\thispagestyle{empty}

\begin{center}
{\bf{\LARGE\raggedright PBW-deformations of smash products involving Hopf algebra of Kac-Paljutkin type}}
\footnotetext {
$\dag$  The corresponding author\\
This work is supported by National Natural Science Foundation of China
(Grant No. 12471038).}


\bigbreak

\normalsize Yujie Gao$^1$,  Shilin Yang$^{{1,2}^{\dag}}$

{\normalsize  \affil[1]{\text{School of Mathematics, Statistics and Mechanics},  
          \text{Beijing University of Technology}, \\
         \text{Beijing}, 
          \text{100124}, 
          \text{P. R. China}}}
          
 {\normalsize   \affil[2]{\text{Department of General Education, Shandong Xiehe University, Jinan},\\
   \text{Shandong}, \text{250109},  \text{P. R. China}}}\\
\end{center}
\begin{quote}
{\noindent{\bf Abstract.}
 Let $H_{2n^2}$ be the Kac-Paljutkin type Hopf algebra of dimension $2n^2$, 
 $A$ its graded Koszul Artin-Schelter regular $H_{2n^2}$-module algebra of  dimension $2$, $A^!$ the Koszul dual of $A$, and $A^{\mathrm{op}}_c$ the braided-opposite algebra of  $A$.  This paper describes  $(0, 1)$-degree PBW-deformations of the smash product $A \sharp H_{2n^2}$ and those of $A^! \sharp\, H_{2n^2}$ under the condition that the Koszul dual $A^!$ of $A$ is also an $H_{2n^2}$-module algebra.
 Also, $0$-degree PBW-deformations of
 $(A \otimes^c A^{\mathrm{op}}_c) \sharp\, H_{2n^2}$ are explored, where $A \otimes^c A^{\mathrm{op}}_c$ is the associated braided tensor product algebra.



}

\end{quote}


\everymath{\displaystyle}


\section*{Introduction}\label{sec:intro}
The PBW-deformation theory for classical Koszul algebras has been extensively studied by numerous authors (see, e.g., \cite{BG, PP, FV}). 
In \cite{CS}, Cassidy and Shelton utilized central extensions to establish the so-called Jacobi condition and provided sufficient and necessary  criteria for PBW-deformations of general connected graded algebras. Subsequently, Xu and Yang provided an algorithm for PBW-deformations of quantum groups based on the Jacobi condition in \cite{XY}.
Recently, a general framework that  includes several versions of the classical PBW-theorem was introduced in \cite{AS}. 
 For more results related to PBW-deformations of some graded algebras, one can refer to 
\cite{HVZ,CEWW,HV, SG, XZ}, among others.
 
Walton and Witherspoon in \cite{WW1} constructed a free resolution of the smash product \(A \sharp\, H\) and established sufficient and necessary conditions for PBW-deformations of \(A \sharp\, H\), where \(H\) is a Hopf algebra and \(A\) is a Koszul \(H\)-module algebra. They applied the obtained results to various Hopf algebras acting on some Koszul algebras, and got corresponding nontrivial PBW-deformations. 
These works generalized the results of \cite{BG, SW}.
Additionally, Walton and Witherspoon in \cite{WW2} described the structures of the braided products \(A \otimes^c B\) as \(H\)-module algebras and provided several examples of PBW-deformations of \((A \otimes^c B) \sharp\, H\).
Their constructions generalized several double constructions described in the literature, including those for Weyl algebras and some classes of Cherednik algebras (see \cite{BB1, BB2}). 

Kac and Paljutkin \cite{KP} discovered  a semisimple Hopf algebra $H_8$ of dimension 8, now named
Kac-Paljutkin Hopf algebra,  which is both noncommutative and noncocommutative.
As a generalization,
Pansera \cite{PD} utilized  Hopf Ore extensions of automorphism type
to construct a class of $2n^2$-dimensional noncommutative and noncocommutative semisimple Hopf algebras $H_{2n^2}$. 
Subsequently, the Grothendieck rings of these algebras were explicitly described using generators and relations in \cite{CYW}. 
Meanwhile, it was demonstrated in \cite{FKMW3} that each Hopf algebra $H_{2n^2}$ serves as a reflection Hopf algebra for an Artin-Schelter regular algebra (or shortly AS regular algebra) of  dimension $2$.

Based on the AS regular $H_{2n^2}$-module algebras $A$ of  dimension $2$  presented 
in \cite{FKMW3}, we aim to describe nontrivial $(0, 1)$-degree PBW-deformations of the smash products
 $A \sharp\, H_{2n^2}$, which extend the discussions for the 8-dimensional Kac-Paljutkin Hopf
 algebra $H_8$ in \cite{WW1}.  Also, we   describe $(0, 1)$-degree PBW-deformations of
  $A^! \sharp\, H_{2n^2}$
under the condition that the Koszul dual $A^!$ is an $H_{2n^2}$-module algebra. Finally, we focus on $(0, 1)$-degree PBW-deformations 
  of the smash products arising from braided product algebras, on  which the Hopf algebra 
  $H_{2n^2}$ acts. These works may provide us with more examples and much information to further understand and develop the theory of PBW-deformations of the smash products in better ways.

Now, we give the outline of the paper.
In Section \ref{sect-1}, we give some basic definitions and results that we need.
For example, the Hopf algebra $H_{2n^2}$ of Kac-Paljutkin type, the PBW-deformation theorem for smash product, and
so on. 
In Section \ref{sec-2},  we describe $(0, 1)$-degree  PBW-deformations of the smash products
based on the Koszul $H_{2n^2}$-module algebras $A^{\pm}_u$ and $B^{\pm}_{u}$ given 
in \cite{FKMW3}.   In Section \ref{sec-4},  
we study Koszul dual algebras of the $H_{2n^2}$-module algebras $ A^{\pm}_u $ and $ B^{\pm}_u $ and calculate $(0, 1)$-degree PBW-deformations of $(A^{\pm}_u)^! \sharp\, H_{2n^2}$ and $(B^{\pm}_u)^! \sharp\, H_{2n^2}$ under the conditions that $(A_u^{\pm})^!$ and $(B^{\pm}_u)^!$
are also $H_{2n^2}$-module algebras,
respectively.
 In Section \ref{sec-3},  
applying 
the braiding $c = J \circ \tau$ for a fixed Drinfeld twisting  $J$ of $H_{2n^2}$,
we determine the braided products $A^{\pm}_u \otimes^c (A^{\pm}_u)_c^{\mathrm{op}}$ and $B^{\pm}_u \otimes^c (B^{\pm}_u)_c^{\mathrm{op}}$, and provide $0$-degree PBW-deformations of 
$\left(A^{\pm}_u \otimes^c (A^{\pm}_u)_c^{\mathrm{op}}\right) \sharp\, H_{2n^2}$ and $\left( B^{\pm}_u \otimes^c (B^{\pm}_u)_c^{\mathrm{op}}\right) \sharp\, H_{2n^2}$.
We emphasize that the techniques used  here are based on the results in \cite{WW1}.

Throughout this paper, all works are over an algebraic closed field $\mathbb{K}$.
For  a Hopf algebra $H$, the Sweedler's notation is used. For example
$$\Delta(h) = \sum_{(h)} h_{(1)} \otimes h_{(2)}$$
for $h\in H$.
The following notations are also used.
\begin{enumerate}
  \item[] $\otimes:$  \,\qquad the tensor product over the field $\mathbb{K}$.
  \item[] $\mathbb{N}:$  \,\qquad the set $\{0, 1, 2, \cdots \}.$
  \item[] $\Z_n$: \,\qquad the residue ring of $\Z$ modulo $n$.
 \item[] $\delta_{a,b}:$ \, the Kronecker sign.
\end{enumerate}

\section{Preliminaries}\label{sect-1}

We always assume that $n>1$ and
$q=e^{\frac{2\pi i}{n}}$, a fixed primitive $n$-th root of unity.
We first recall  the Hopf algebra $H_{2n^2}$, which is a generalization of the Kac-Paljutkin Hopf algebra $H_8.$

\begin{defn}\cs{PD} The  Hopf algebra $H_{2n^2}$ is  the associative algebra generated by $x,y$ and $z$, with the following relations
\begin{equation*}\label{eqn1-1}
   \begin{array}{lll}
    &x^n=1,\quad y^n=1 \\
     &xy=yx,\quad zx=yz,\quad zy=xz, \\
    &z^2=\frac{1}{n}\sum\limits_{i,j=0}^{n-1}q^{-ij}x^iy^j.
   \end{array}
\end{equation*}
The coproduct $\Delta$ of $H_{2n^2}$ is given by
$$\Delta(x)=x\otimes x,\quad \Delta(y)=y\otimes y,\quad \Delta(z)=\frac{1}{n}\sum\limits_{i,j=0}^{n-1}q^{-ij}x^{i}z \otimes y^j z,$$
 the counit $\epsilon$ is given by
$$\epsilon(x)=1, \quad \epsilon(y)=1,\quad \epsilon(z)=1,$$
and the antipode $S$ is  given by
$$S(x)=x^{-1}, \quad S(y)=y^{-1},\quad \quad S(z)=z. $$
\end{defn}
The Hopf algebra $H_{2n^2}$ is  semisimple and has dimension $2n^2$ with a basis
  $$\{x^iy^j,x^iy^jz| 0 \leq i,j\leq n-1\}.$$
Here we  list irreducible modules  of $H_{2n^2}$ according to \cite{ CYW, FKMW3}.
For the purpose, we introduce the notation
$p=-e^{\frac{\pi i}{n}} = e^{\frac{(n+1)\pi i}{n}}$. 
Notably, when $n$ is odd, $p$ is a primitive $n$-th root of unity.
\begin{prop} {\rm (see \cite{ CYW, FKMW3})} \label{prop1-2}
The set
$$\mathcal{S}=\left\{S_k^{\pm}| k \in \mathbb{Z}_n\right\}\cup \left\{S_{k,l} | 0 \leq k < l \leq n-1\right\}$$
form a complete set of irreducible modules of the Hopf algebra $H_{2n^2}$, where  $S_k^{\pm} \ (k\in \mathbb{Z}_n)$ is the one-dimensional irreducible $H_{2n^2}$-module with a basis $\{u^k\}$, on which the actions are
$$x\cdot u^k=q^k u^k,\quad y\cdot u^k=q^k u^k,\quad z\cdot u^k=\pm p^{k^2} u^k;$$
$S_{k,l} (0\leq k<l \leq n-1)$ is the two-dimensional irreducible $H_{2n^2}$-module
with a basis $\{u_1^{kl},\ u_2^{kl}\}$, on which the actions are
\begin{alignat*}{3}
    &x\cdot u_1^{kl}=q^k u_1^{kl},
    \quad\quad&&y\cdot u_1^{kl}=q^lu_1^{kl},
    \quad\quad&&z\cdot u_1^{kl}
    =q^{kl}u_2^{kl},\\
    &x\cdot u_2^{kl}=q^l u_2^{kl},
    \quad\quad&&y\cdot u_2^{kl}=q^k u_2^{kl},
    \quad\quad&&z\cdot u_2^{kl}
    =u_1^{kl}.
    \end{alignat*}
\end{prop}

Let $V$ be a finite-dimensional vector space, and
\[T(V) = \mathbb{K} \oplus V \oplus V^{\otimes 2} \oplus \cdots\]
be the corresponding connected tensor algebra. An $\mathbb{N}$-graded $\mathbb{K}$-algebra $A = T(V) / (I)$ is called a quadratic algebra if it is generated by  elements in $V$ of degree one 
with  the quadratic relations $I \subseteq V \otimes V$. 
Here, we identify $I$ with the subspace $\sum_{r_i\in I} \mathbb{K}r_i$ of $V\otimes V.$

The connected graded algebra $A$ is called Koszul if the trivial module $_A\mathbb{K}$ admits a resolution

\[
\cdots \rightarrow P_i \rightarrow \cdots \rightarrow P_1 \rightarrow P_0 \rightarrow {_A\mathbb{K}} \rightarrow 0
\]
such that each $P_i$ is a  graded free left $A$-module generated in degree $i$.

 A connected graded algebra $A$ is said to be an Artin-Schelter regular algebra (shortly AS regular algebra)  of dimension $d$
 if $A$ has finite global dimension $d$, finite GK dimension, and is AS-Gorenstein.

Now, let $H$ be a Hopf algebra, and $V$ an  left $H$-module.
It is well-known that $H$, $V\otm H$  are $H$-modules defined by 
$$h\cdot_{\rm ad} l=\sum_{(h)} h_{(1)}lS(h_{(2)}), \hbox{ (the left adjoint action),}$$
and
$$h\cdot_{\rm ad}( v\otm l)=\sum_{(h)} h_{(1)}\cdot v\otm h_{(2)}\ad l,$$
for all $h\in H$ and $l\in H,$ and $v\in V$, respectively.

We say that an algebra $A$ is an  $H$-module algebra if
\begin{enumerate}
  \item[(1)] $A$ is a  $H$-module;
  \item[(2)]  $h \cdot (ab) = \sum\limits_{(h)} (h_{(1)} \cdot a)(h_{(2)} \cdot b) \quad \text{and} \quad h \cdot 1_A = \epsilon(h) 1_A $
for all $ h \in H $, $ a, b \in A $.
\end{enumerate}
 
 Let $H$ be a Hopf algebra and $A$ an $H$-module algebra, the smash product $A\sharp\,H$ is 
 the associative algebra defined by
 \begin{enumerate}
   \item[(1)] as a $\K$-space, $A\sharp\,H=A\otimes H$;
   \item[(2)] the multiplication is
   $$(a\sharp\,g)(b\sharp\,h)=\sum_{(g)} a(g_{(1)}\cdot b)\sharp\,g_{(2)}h,$$
   for all $a, b\in A$ and $g, h\in H$.
 \end{enumerate}

Suppose $A = T(V)/(I)$ is an $\mathbb{N}$-graded, Koszul, (left) $H$-module algebra, where $I \subseteq V \otimes V$, and the actions of the Hopf algebra $H$ preserve the grading of $A$.
 It is straightforward to see that the subspace $V$ is an $H$-submodule of $A$. Additionally, the action of $H$ on $V$ can be naturally extended via the coproduct of $H$ such that
 $T(V)$ is an $H$-module algebra. 
Since $I$ is precisely the subspace of degree 2 in $T(V)$, it must be preserved under the action of $H$. Therefore, $I$ is also an $H$-submodule of $T(V)$.

Let $\mathcal{C}$ denote a monoidal category of (left) $H$-modules  for a Hopf algebra $H.$ 
We say that  $\mathcal{C}$ is equipped with a braiding $\mathrm{c}$, if there exist functorial isomorphisms
$$
\mathrm{c}_{A, B}: A \otimes B \stackrel{\sim}{\longrightarrow} B \otimes A
$$
for all pairs of objects $A, B$ in $\mathcal{C}$, satisfying the hexagon axioms (see for example \cite{KC}).

\begin{defn} {\rm (see \cite{WW2})} \label{defn1-3}
 Let $A, \ A'$ be two (left) $H$-module algebras in the monoidal category $\mathcal{C}$.
 \begin{itemize}
   \item [(a)] The braided product $A \otimes^{\mathrm{c}} A'$ is
the object  $A\otimes A'$ of $\mathcal{C}$ as a $\mathbb{K}$-space,  which multiplication is defined by the composition
$$
A \otimes A' \otimes A \otimes A' \xrightarrow{1 \otimes \mathrm{c} \otimes 1} A \otimes A \otimes A' \otimes A' \xrightarrow{m_A \otimes m_{A'}} A \otimes A'
$$
where $\mathrm{c}=\mathrm{c}_{A, A'}$ and $m_A$ and $m_{A'}$ are multiplication on $A$ and $A'$, respectively.
   \item [(b)] The braided-opposite algebra of the (left) $H$-module algebra $A$, denoted by $A_{\mathrm{c}}^{\mathrm{op}}$, is
       the  object $A$ of $\mathcal{C}$ as a $\mathbb{K}$-space, which multiplication is defined by
$$
A \otimes A \xrightarrow{\mathrm{c}} A \otimes A \xrightarrow{m_A} A
$$
where $\mathrm{c}=\mathrm{c}_{A, A}$.
 \end{itemize}
\end{defn}
The multiplications in the above definition are indeed associative, see \cite{B} and \cite{MS}.

Let
$$ A =\bigoplus_{j\geq 0} A_j, \ A_0 = \mathbb{K} $$
 be a connected graded algebra  given by generators $a_1, \ldots, a_n$ together with some homogeneous relations $r_1, \ldots, r_m.$ Recall that a deformation of  $A$ is an algebra $U$ given by generators $a_1, \ldots, a_n$ with the relations $r_1+ l_1, \ldots, r_m+l_m$, 
 where each $l_j$ is a possibly nonhomogeneous element of degree less than the degree of $r_j$. It is obvious that $U$ is a filtered algebra.
 A PBW-deformation of $A$ is defined to be
  $U$ with an ascending filtration
 $$
 \mathcal{F}_{-1}U=0\subset\mathcal{F}_0 U \subset \mathcal{F}_1U \subset\mathcal{F}_2U \subset \cdots
 $$
such that the associated graded algebra 
${\rm gr}(U)=\oplus_{i\ge 0} \mathcal{F}_i U/ \mathcal{F}_{i-1}U$ 
is isomorphic to the algebra $A$.

\begin{exam} Let $A=\mathbb{K}[x, y]$, then $U=\mathbb{K}\lr{x, y}/(xy-yx-1)$
is a PBW-deformation of $A$.
\end{exam}

\begin{defn}[see {\cite[\bf Notation 0.3]{WW1}}] \label{defn1-4} Let $A=T(V) /(I)$ be an $\mathbb{N}$-graded, Koszul, (left) $H$-module algebra and quadratic relations $I \subseteq V \otimes V$, which the action of $H$ preserves the grading and the subspace $I$ of $V \otimes V$. 
Let $\mathscr{U}_{A, \kappa}$ be the filtered $\mathbb{K}$-algebra  of $A\sharp H$ defined by
$$\mathscr{U}_{A, \kappa}=\frac{T(V) \sharp\, H}{\left( r-\kappa(r)\,|\, {r \in I}\right)},$$
where $\kappa : I\rightarrow H\oplus (V\otimes H)$ is a $\mathbb{K}$-bilinear map, which is decomposed into its constant part $\kappa^{C} : I\rightarrow H$ and linear part $\kappa^{L} : I\rightarrow V\otimes H$.
\end{defn}

By {\rm\cite{WW1}}, elements in $H$ are assigned degree $0$, and elements in $V$ 
degree $1$,
and $A\sharp H\cong (T(V)\sharp H)/{(I)}$. Therefore, $\mathscr{U}_{A, \kappa}$ can  be viewed as
a deformation of $ A\sharp H$ naturally. Accordingly,  to get a deformation $\mathscr{U}_{A, \kappa}$ of $A\sharp H$, which is a filtered algebra, the 
additional parts $\kappa(I)$ should be contained in 
$H\oplus (V\otimes H).$ Thus we can say that $\mathscr{U}_{A, \kappa}$ is  a $(0,1)$-degree deformation of $A\sharp H,$  in particular, a $0$-degree deformation if $\kappa^L=0$.

The following result is known.
\begin{thm}\cs{WW1, BG, SW} \label{thm1-4}
The algebra $\mathscr{U}_{A,\kappa}$ is a $\mathrm{PBW}$-deformation-deformation of $A\sharp\, H$ if and only if the following conditions hold.
\begin{enumerate}
  \item [(1)]$\kappa$ is $H$-invariant\,$:$  $\kappa(h\cdot r)=h\ad\kappa(r)$ for all $r\in I, h\in H;$ 
  \end{enumerate}
and also if $\mathrm{{dim}} V\geq 3$,  the maps $\kappa^C \otimes \mathrm{id}-\mathrm{id} \otimes \kappa^C$ and $\kappa^L \otimes \mathrm{id}-\mathrm{id} \otimes \kappa^L$, which are defined on $(I \otimes V) \cap(V \otimes I),$ satisfy
 \begin{enumerate}
  \item [(2)]$\operatorname{Im}\left(\kappa^L \otimes \mathrm{id}-\mathrm{id} \otimes \kappa^L\right) \subseteq I;$\\
\item [(3)] $\kappa^L \circ\left(\kappa^L \otimes \mathrm{id}-\mathrm{id} \otimes \kappa^L\right)=-\left(\kappa^C \otimes \mathrm{id}-\mathrm{id} \otimes \kappa^C\right);$\\
\item [(4)] $\kappa^C \circ\left(\mathrm{id} \otimes \kappa^L-\kappa^L \otimes \mathrm{id}\right) \equiv 0.$
\end{enumerate}
\end{thm}

Thereafter, in calculations we often omit the tensors $\otm$ that the meaning is 
obvious and without confusions.

Based on the PBW-deformation theory of a quadratic algebra $A$, Theorem  \ref{thm1-4} provides  us a constructive method to  get the $(0, 1)$-degree PBW-deformations of $A\sharp\, H$.

In the rest of the paper, we shall apply Theorem \ref{thm1-4} to explore nontrivial $(0,1)$-degree PBW-deformations of the smash products $A \sharp H_{2n^2}$ and $A^! \sharp H_{2n^2}$, where $A$ is
the graded  Koszul AS regular $H_{2n^2}$-module algebra of  dimension $2$, and $A^!$ is the Koszul dual  of $A$. We also study the nontrivial $0$-degree PBW-deformations of braided tensor product algebra $A \otimes^c A^{\mathrm{op}}_c$.

Here we point out that 
Theorem \ref{thm1-4} can not provide us a method to classify all PBW-deformations up to isomorphism.
To the best of our knowledge, there are few works to classify all such possible deformations.
 Theorem \ref{thm1-4} can be used to construct nontrivial $(0,1)$-degree PBW-deformations.

\section{The PBW-deformations
associated to AS regular algebras}\label{sec-2}
In this section,  we restrict our attention to  the AS regular algebras  of 
dimension $2$, such that 
they are Koszul $H_{2n^2}$-module algebras. 

Recall that the Hopf algebra $H_{2n^2}$ is semisimple with
the complete set $\mathcal{S}$ of irreducible modules of dimensions $1$ and $2$ (see Proposition \ref{prop1-2}). Based on this fact,  there only exist two classes of  AS regular (left) 
$H_{2n^2}$-module algebras of  dimension $2$ by  \cite{FKMW3} : 
\begin{eqnarray*}
 A_u^{\pm, k, l} &=& \frac{\mathbb{K}\langle u_1, u_2 \rangle}{(u_1 u_2 \pm p^{l^2 - k^2} u_2 u_1)},  \hbox{ where } 0 \leq k < l \leq n-1;\\
B_{u}^{\pm, k, l}&=&\frac{\mathbb{K}\langle u_1,u_2\rangle}{\left( u_1^2\pm q^{kl}u_2^2\right)}, \hbox{ where } n \hbox{ is even, }  0\leq k<l\leq n-1, \hbox{ and } \ k\equiv l\ (\mathrm{mod}\ \textstyle\frac{n}{2}),
\end{eqnarray*}
where the  module structures of $A_u^{\pm, k, l}$ ( resp, 
$B_{u}^{\pm, k, l}$) are given by
\begin{equation}\label{action}
\begin{aligned}
    &x\cdot u_1=q^k u_1,
  \quad&&y\cdot u_1=q^lu_1,
   \quad&&z\cdot u_1
    =q^{kl}u_2,\\
    &x\cdot u_2=q^l u_2,
   \quad&&y\cdot u_2=q^k u_2,
    \quad&&z\cdot u_2
    =u_1,
    \end{aligned}
    \end{equation}
deduced from $S_{k,l}$
in Proposition 1.2

In this paper, we only explore those of $A_u^{\pm, k, l}$ ( resp. 
$B_{u}^{\pm, k, l}$) for some suitable fixed numbers $k$ and $l$. So, by abuse of notations, we omit $k, l$ in the superscript, that is $A_u^{\pm}:=A_u^{\pm, k, l}\hbox{ and } B_{u}^{\pm}:=B_{u}^{\pm, k, l}.$
By \cite[Sect. 2]{C}, the AS regular algebra of  dimension $2$ generated by elements of degree $1$ is, up to isomorphism, either a quantum plane or a Jordan plane, both of which are known to be Koszul. 
Therefore, $A_u^{\pm}$ are $\mathbb{N}$-graded Koszul 
$H_{2n^2}$-module algebras. 
We also claim that $B_{u}^{\pm}$ are   $\mathbb{N}$-graded Koszul  $H_{2n^2}$-module algebras.
To see this,  it is straightforward to verify that 
the algebra $B_{u}^{+}$ (resp. $B_{u}^{-}$) admits an isomorphism to the algebra $\frac{\mathbb{K}\langle u_1,u_2\rangle}{\left( u_1u_2+u_2u_1\right)}$, which is  an $\mathbb{N}$-graded Koszul algebra by  \cite[Example 4.2]{WW1}. These isomorphisms $f_1$ (resp, $f_2$) are  as follows:
\begin{align*}
   &f_1:\quad u_1 \mapsto {\bf i}p^{kl}u_1-\frac{{\bf i}
   p^{-kl}}{2}u_2 , \quad u_2 \mapsto u_1+\frac{q^{-kl}}{2}u_2,\\
 & f_2:\quad u_1 \mapsto p^{kl}u_1+\frac{p^{-kl}}{2}u_2 , \quad u_2 \mapsto u_1-\frac{q^{-kl}}{2}u_2,
\end{align*}
where $\mathbf{i}=\sqrt{-1}.$ 
However, the algebras $B_u^\pm$ are not isomorphic to $\frac{\mathbb{K}\langle u_1,u_2\rangle}{\left( u_1u_2+u_2u_1\right)}$ as $H_{2n^2}$-module algebras.
We fix the above notations in subsequent parts.

Recall that the left $H_{2n^2}$-adjoint action on $H_{2n^2}$ is
given by
$$h\ad \ell:=\sum\limits_{(h)} h_{(1)}\ell S(h_{(2)}),$$
for all $h,\ell\in H_{2n^2}$.  

\begin{lem}\label{lem2-3} The adjoint action of $z$ on $\sum\limits_{i,j=0}^{n-1}a_{ij}x^iy^j,$ where $a_{ij}\in \mathbb{K}$
is
$$z\ad\left(\sum\limits_{i,j=0}^{n-1}a_{ij}x^iy^j\right)=\sum\limits_{i,j=0}^{n-1}a_{ji}x^iy^j.$$
\end{lem}
\begin{proof} It is straightforward to see that
\begin{align*}
 z\ad\left(\sum\limits_{i,j=0}^{n-1}a_{ij}x^iy^j\right)=
  &\frac{1}{n}\sum\limits_{g,h=0}^{n-1} q^{-gh}x^gz\left(\sum\limits_{i,j=0}^{n-1}a_{ij}x^iy^j\right)zy^{-h} \\
  = &\frac{1}{n}\sum\limits_{g,h=0}^{n-1} q^{-gh} \sum\limits_{i,j=0}^{n-1}a_{ij}x^{g+j}y^{i-h}z^2\\
 =&\frac{1}{n}\sum\limits_{v,w=0}^{n-1}  \sum\limits_{i,j=0}^{n-1}a_{ij}q^{-(v-j)(i-w)}x^{v}y^{w}z^2\\
 =&\frac{1}{n^2}\sum\limits_{v,w=0}^{n-1}  \sum\limits_{i,j=0}^{n-1}\sum\limits_{r,m=0}^{n-1}\left(a_{ij}q^{-(v-j)(i-w)-rm}\right)x^{v+r}y^{w+m}\\
= &\frac{1}{n^2}\sum\limits_{s,t=0}^{n-1}\left( \sum\limits_{v,w=0}^{n-1}  \sum\limits_{i,j=0}^{n-1}a_{ij}q^{-(v-j)(i-w)-(s-v)(t-w)}\right)x^{s}y^{t}\\
= &\frac{1}{n^2}\sum\limits_{s,t=0}^{n-1}q^{-st}\left(\sum\limits_{i,j=0}^{n-1}a_{ij}q^{ij}  \sum\limits_{w=0}^{n-1}q^{w(s-j)}\sum\limits_{v=0}^{n-1}q^{v(t-i)}\right)x^{s}y^{t}\\
=&\sum\limits_{i,j=0}^{n-1}a_{ji}x^iy^j.
\end{align*}
The proof is completed.
\end{proof}
\begin{lem}\label{lem2-4}
The adjoint action of $z$ on $\sum\limits_{i,j=0}^{n-1}b_{ij}x^iy^jz,$ where $b_{ij}\in \mathbb{K}$
is
$$z\ad\left(\sum\limits_{i,j=0}^{n-1}b_{ij}x^iy^jz\right)= \sum\limits_{s,t=0}^{n-1}\left(\frac{1}{n}\sum\limits_{i,j=0}^{n-1}b_{ij}q^{(t-i)(i+j-s-t)}  \right)x^{s}y^{t}z.$$
\end{lem}
\begin{proof}  It is straightforward to see that
\begin{align*}
 z\ad\left(\sum\limits_{i,j=0}^{n-1}b_{ij}x^iy^jz\right)=
  &\frac{1}{n}\sum\limits_{g,h=0}^{n-1} q^{-gh}x^gz\left(\sum\limits_{i,j=0}^{n-1}b_{ij}x^iy^jz\right)zy^{-h} \\
  = &\frac{1}{n}\sum\limits_{g,h=0}^{n-1} q^{-gh} \sum\limits_{i,j=0}^{n-1}b_{ij}x^{g+j-h}y^{i}z^3\\
 =&\frac{1}{n}\sum\limits_{v,h=0}^{n-1}  \sum\limits_{i,j=0}^{n-1}b_{ij}q^{-(v-j+h)h}x^{v}y^{i}z^3\\
 =&\frac{1}{n^2}\sum\limits_{v,h=0}^{n-1}  \sum\limits_{i,j=0}^{n-1}\sum\limits_{r,m=0}^{n-1}b_{ij}q^{-(v-j+h)h-rm}x^{v+r}y^{i+m}z\\
= &\frac{1}{n^2}\sum\limits_{s,t=0}^{n-1}\left( \sum\limits_{v,h=0}^{n-1}  \sum\limits_{i,j=0}^{n-1}b_{ij}q^{-(v-j+h)h-(s-v)(t-i)}\right)x^{s}y^{t}z\\
= &\frac{1}{n^2}\sum\limits_{s,t=0}^{n-1}\left(\sum\limits_{i,j=0}^{n-1}b_{ij}q^{-s(t-i)}  \sum\limits_{h=0}^{n-1}q^{h(j-h)}\sum\limits_{v=0}^{n-1}q^{v(t-i-h)}\right)x^{s}y^{t}z\\
=&\sum\limits_{s,t=0}^{n-1}\left(\frac{1}{n}\sum\limits_{i,j=0}^{n-1}b_{ij}q^{(t-i)(i+j-s-t)}  \right)x^{s}y^{t}z.
\end{align*}
The proof is completed.
\end{proof}
For the $\mathbb{N}$-graded Koszul $H_{2n^2}$-module algebras $A_u^\pm=T(V)/(I)$, where $V=\mathbb{K}u_1+\mathbb{K}u_2$ and $$I=\{r:=u_1u_2\pm p^{l^2-k^2}u_2u_1\}.$$
According to (\ref{action}), we have
$$x\cdot r=q^{k+l}r,\quad y\cdot r=q^{k+l}r,\quad z\cdot r=\pm p^{(k+l)^2}r,$$
where  $0 \leq k< l \leq n-1.$

The smash products
$A_{u}^{\pm}\sharp\, H_{2n^2}$ can be viewed as the algebras generated by $x, y, z, u_1, u_2$ with the relations
\begin{eqnarray*}\label{eqn1-1}
   \begin{array}{lll}
    &x^n=1,\quad y^n=1, \\
     &xy=yx,\quad zx=yz,\quad zy=xz, \\
    &z^2=\frac{1}{n}\sum\limits_{i,j=0}^{n-1}q^{-ij}x^iy^j,\quad u_1u_2=\pm p^{l^2-k^2}u_2u_1,\\
&xu_1=q^ku_1x, \quad xu_2=q^lu_2x,\\
&yu_1=q^lu_1y, \quad yu_2=q^ku_2y,\\
&zu_1=q^{kl}u_2y^{l}z, \quad zu_2=u_1y^{k}z,\\
   \end{array}
\end{eqnarray*}
where  $0 \leq k< l \leq n-1$.

Let $\kappa : I\rightarrow H_{2n^2}\oplus (V\otimes H_{2n^2})$ be a $\mathbb{K}$-bilinear map as stated in Definition \ref{defn1-4}, where $\kappa$ is the sum of its constant part $\kappa^{C}$ and its 
linear part $\kappa^{L}$.
  We get the first main result  as follows.
\begin{thm}\label{thm2-5}Keep notations as above. Let $A_u^\pm=T(V)/(I)$, where $V=\mathbb{K}u_1+\mathbb{K}u_2$ and $I=\{r:=u_1u_2\pm p^{l^2-k^2}u_2u_1\}.$
For  $0 \leq k< l \leq n-1$ ,  the deformation
$$\mathscr{U}_{A_{u}^{\pm}, \ \kappa} = \frac{\mathbb{K}\langle u_1,u_2 \rangle  \sharp\, H_{2n^2}}{\left( r-\kappa(r)\,|\, {r \in I}\right)},$$
is a   nontrivial  $\mathrm{PBW}$-deformation of $A_{u}^{\pm}\sharp\, H_{2n^2}$
if and only if  one of the following holds$:$
\begin{enumerate}
  \item [{(a)}] $k+l\equiv 0\ (\mathrm{mod}\ n)$, and
  $$\begin{aligned}
&\kappa(r)=\sum_{i=0}^{n-1}{\mathbf\delta}_{+,\pm}\, a_{ii}x^iy^i
+ \sum_{0 \leq i<j \leq n-1} a_{ij}(x^iy^j \pm x^jy^i)
+ \sum_{i,j=0}^{n-1}\delta_{+,\pm} b_{i+j} x^iy^jz\\
&+\sum\limits_{i,j=0}^{n-1}q^{ik}c_{i+j} u_1 x^iy^jz
    \mathbf{\pm}  \sum\limits_{i,j=0}^{n-1}q^{jk}{c_{i+j+k}} u_2x^iy^jz;\\
\end{aligned}$$
  \item [{(b)}] $n\ \text{ is even},\ k+l\not\equiv 0\ (\mathrm{mod}\ n), \ k+l\equiv 0\ (\mathrm{mod}\ \textstyle\frac{n}{2})$, and
  $$
  \kappa(r)=\pm\sum_{i,j=0}^{n-1}(-1)^jp^{-\frac{n^2}{4}}b_{i+j-
\frac{n}{2}}x^iy^jz,
$$
\end{enumerate}
  where $a_{ij}, b_i, c_i \in \mathbb{K}$ for all $i, j \in \mathbb{Z}_n$, not all of which are zero.
\end{thm}

\begin{proof} To get the PBW-deformations of $A_{u}^{\pm}\sharp\, H_{2n^2}$, it is sufficient to consider the condition of Theorem (\ref{thm1-4})(1)
since $\dim V = 2$. Recall that $\kappa=\kappa^C+\kappa^L$, we
now begin by computing $\kappa^C$. 

Let 
 $$\kappa^C(r)=\sum\limits_{i,j=0}^{n-1}\left(a_{ij}x^iy^j+b_{ij}x^iy^jz\right),$$
where $a_{ij},b_{ij}\in \mathbb{K}$.

Since
$$h\ad(\kappa^C(r))=\sum_{(h)} h_{(1)}\kappa^C(r)S(h_{(2)}),$$
 for all $h\in H_{2n^2}$, it yields  that
\begin{align}\label{align2-1}
  x\ad \kappa^C(r)
             =&\sum\limits_{i,j=0}^{n-1}\left(a_{ij}x^{i}y^j+b_{ij}x^{i+1}y^{j-1}z\right),
\end{align}
\begin{align}\label{align2-2}
 y\ad \kappa^C(r)=\sum\limits_{i,j=0}^{n-1}\left(a_{ij}x^{i}y^j+b_{ij}x^{i-1}y^{j+1}z\right).
 \end{align}
The conditions that $x\ad \kappa^C(r)=\kappa^C(x\cdot r)=q^{k+l}\kappa^C(r)$ and  $y\ad \kappa^C(r)=\kappa^C(y\cdot r)=q^{k+l}\kappa^C(r)$ imply that
\begin{equation}\label{eqn2-4-1}
a_{ij}=q^{k+l}a_{ij},\quad b_{i-1,j+1}=q^{k+l}b_{i,j} \quad \text{and}\quad  b_{i+1,j-1}=q^{k+l}b_{i,j},
\end{equation}
for all $i, j\in\Z_n$.

Let us consider the matrix $A=(a_{ij})_{n\times n}$ (resp. $B=(b_{ij})_{n\times n}$), where $i, j\in\Z_n.$

{\bf (a).  $k+l\equiv 0\ (\mathrm{mod}\ n)$.}

 We observe that
\begin{equation*}\label{eqn2-1}
    b_{i+1,j-1} = b_{i,j} = b_{i-1,j+1}, \quad   \text{for\ all}\ i, j\in\Z_n,
\end{equation*}
and 
\begin{equation}\label{eqn2-2}
b_{ij}=b_{0, i+j},  \quad \text{for \ all} \  i, j\in\Z_n.
\end{equation}
Since
 \begin{align*}
 z\ad(\kappa^C(r))
 =&\sum\limits_{i,j=0}^{n-1}a_{ji}x^iy^j+\sum\limits_{s,t=0}^{n-1}\left(\frac{1}{n}\sum\limits_{i,j=0}^{n-1}b_{ij}q^{(t-i)(i+j-s-t)}\right)x^sy^tz\\ =&\sum\limits_{i,j=0}^{n-1}a_{ji}x^iy^j+\sum\limits_{s,t=0}^{n-1}\left(\frac{1}{n}\sum\limits_{i,j=0}^{n-1}b_{0,i+j}q^{t(i+j-s-t)}q^{i(s+t-i-j)}\right)x^sy^tz\\ =&\sum\limits_{i,j=0}^{n-1}a_{ji}x^iy^j+\sum\limits_{s,t=0}^{n-1}\left(\frac{1}{n}\sum\limits_{i,r=0}^{n-1}b_{0,r}q^{t(r-s-t)}q^{i(s+t-r)}\right)x^sy^tz\\  =&\sum\limits_{i,j=0}^{n-1}a_{ji}x^iy^j+\sum\limits_{s,t=0}^{n-1}\left(\frac{1}{n}\sum\limits_{r=0}^{n-1}b_{0,r}q^{t(r-s-t)}\sum\limits_{i=0}^{n-1}q^{i(s+t-r)}\right)x^sy^tz\\
 =&\sum\limits_{i,j=0}^{n-1}a_{ji}x^iy^j+\sum\limits_{s,t=0}^{n-1}b_{0,s+t}x^sy^tz,
 \end{align*}
 by Lemma \ref{lem2-3},  Lemma \ref{lem2-4} and 
 $(\ref{eqn2-2})$.

The condition that $z\ad\kappa^C(r)=\kappa^C(z\cdot r)=\kappa^C(r)$ implies that
\begin{equation*}\label{eqn2-3}
 a_{ij}= a_{ji},
\end{equation*}
for all $i, j\in\Z_n$, which means that the matrix $A=(a_{ij})$ is  a symmetric matrix
and
\begin{equation*}
B={\left( {\begin{array}{*{20}{c}}
    b_{0,0} & b_{0,1} & b_{0,2} & \ldots & b_{0,n-2} & b_{0,n-1} \\
    b_{0,1} & b_{0,2} & b_{0,3} & \ldots & b_{0,n-1} & b_{0,0} \\
    b_{0,2} & b_{0,3} & b_{0,4} & \ldots & b_{0,0} & b_{0,1} \\
    \vdots & \vdots & \vdots & \ddots & \vdots & \vdots \\
    b_{0,n-2} & b_{0,n-1} & b_{0,0} & \ldots & b_{0,n-4} & b_{0,n-3} \\
    b_{0,n-1} & b_{0,0} & b_{0,1} & \ldots & b_{0,n-3} & b_{0,n-2} \\
\end{array}} \right).}
\end{equation*}
Consequently,

\begin{equation*}\label{r1}
\kappa^C(r)=\sum\limits_{0 \leq i\leq n-1}a_{ii}x^iy^i+\sum\limits_{0 \leq i<j \leq n-1}a_{ij}(x^iy^j+x^jy^i)+\sum\limits_{0 \leq i,j \leq n-1}b_{0,i+j} x^iy^jz.\end{equation*}

The condition that $z \ad \kappa^C(r) =\kappa^C(z\cdot r)= -\kappa^C(r)$ implies that
\begin{equation*}
a_{ij} = -a_{ji}, \quad b_{ij} = 0, \quad \text{for\ all }\ i, j\in\Z_n,
\end{equation*}
which means that the matrix $A = (a_{ij})$ is an antisymmetric matrix and $B = (b_{ij})_{n \times n} = 0$.
Consequently,
\begin{equation*}\label{2r}
\kappa^C(r)=\sum\limits_{0 \leq i<j \leq n-1}a_{ij}(x^iy^j-x^jy^i).
\end{equation*}

{\bf (b).  $k+l\not \equiv 0\ (\mathrm{mod}\ n)$}.

In this case, the matrix $A = (a_{ij})_{n \times n} = 0$, and
$q^{2(k+l)}b_{ij} =b_{ij}$  by (\ref{eqn2-4-1}). And it is easily seen that 
$\kappa^L(r)=0$ (see the part $\kappa^L(r)$ below).
To obtain nontrivial PBW-deformations of $A_u^\pm \sharp\, H_{2n^2}$,
 we must have $q^{2(k+l)}=1,$ 
which implies that $n$ has to be an even number and $k + l \equiv 0 \pmod{\textstyle\frac{n}{2}}$. Thus, we get  the following relations
\begin{equation*}\label{eqn2-5}
b_{i+1,j-1}=-b_{i,j}= b_{i-1,j+1}, \quad   \text{for\ all } \ i, j\in\Z_n.
\end{equation*}
From these relations, we get
\begin{equation}\label{eqn2-6}
b_{ij}=(-1)^{i}b_{0, i+j}, \quad \text{for\ all}  \ i, j\in\Z_n,
\end{equation}
and
\begin{equation*}B=
{\left( {\begin{array}{*{20}{c}}
    b_{0,0} & b_{0,1} & b_{0,2} & \ldots & b_{0,n-2} & b_{0,n-1} \\
    -b_{0,1} & -b_{0,2} & -b_{0,3} & \ldots & -b_{0,n-1} & -b_{0,0} \\
    b_{0,2} & b_{0,3} & b_{0,4} & \ldots & b_{0,0} & b_{0,1} \\
    \vdots & \vdots & \vdots & \ddots & \vdots & \vdots \\
    b_{0,n-2} & b_{0,n-1} & b_{0,0} & \ldots & b_{0,n-4} & b_{0,n-3} \\
    -b_{0,n-1} & -b_{0,0} & -b_{0,1} & \ldots & -b_{0,n-3} & -b_{0,n-2} \\
\end{array}} \right).}
\end{equation*}
By Lemma \ref{lem2-4} and (\ref{eqn2-6}), we also have
 \begin{align*}
 z\ad(\kappa^C(r))
 =&\sum\limits_{i,j=0}^{n-1}(-1)^jb_{0,i+j-\frac{n}{2}}x^iy^jz.
 \end{align*}
Therefore, the condition that $z\ad \kappa^C(r)=\kappa^C(z\cdot r)=\pm p^{(k+l)^2}\kappa^C(r)$ implies that
\[b_{0,i+j}=\pm (-1)^{i+j} p^{-\frac{n^2}{4}}b_{i+j-\textstyle\frac{n}{2}}.\]
Hence $$\kappa^C(r)=\pm\sum\limits_{i,j=0}^{n-1}(-1)^jp^{-\frac{n^2}{4}}
b_{0,i+j-\frac{n}{2}}x^iy^jz.$$

Subsequently, we try to compute the part $\kappa^{L}(r).$ 

We see that
$\kappa^{L}(r)=u_1g+u_2g' \in V \otimes H$, and

  $$g=\sum\limits_{i,j=0}^{n-1}(d_{ij} x^iy^j+c_{ij} x^iy^jz),  \quad
g'=\sum\limits_{i,j=0}^{n-1}(e_{ij} x^iy^j+f_{ij} x^iy^jz).$$
   Here, $c_{ij}, d_{ij}, e_{ij}, f_{ij}\in \mathbb{K}.$
Assume that $h\in H_{2n^2},$  then
   $$h\ad (\kappa^{L}(r))=\sum_{(h)} (h_{(1)}\cdot u_1)(h_{(2)}gS(h_{(3)})+\sum_{(h)} (h_{(1)}\cdot u_2)\left(h_{(2)}g'S(h_{(3)}\right)$$
by the definition of the left $H_{2n^2}$-adjoint action on
$V \otimes H.$
So we obtain  that
\begin{align*}
x\ad \kappa^{L}(r)
                    = &\sum\limits_{i,j=0}^{n-1}q^{k}u_1(d_{ij}x^iy^j +c_{ij} x^{i+1}y^{j-1}z)
+\sum\limits_{i,j=0}^{n-1}q^{l}u_2(e_{ij} x^iy^j+f_{ij} x^{i+1}y^{j-1}z),
\end{align*}
\begin{align*}
y\ad \kappa^{L}(r)
                    = &\sum\limits_{i,j=0}^{n-1}q^{l}u_1(d_{ij} x^iy^j+c_{ij}x^{i-1}y^{j+1}z)
+\sum\limits_{i,j=0}^{n-1}q^{k}u_2(e_{ij} x^iy^j+f_{ij} x^{i-1}y^{j+1}z),
\end{align*}
by (\ref{align2-1}) and (\ref{align2-2}).
The conditions that $x\ad \kappa^L(r)=\kappa^L(x\cdot r)$ and $y\ad \kappa^L(r)=\kappa^L(y\cdot r)$  imply that
\begin{equation*}
\begin{aligned}
  &d_{ij}=q^{l} d_{ij},\quad d_{ij}=q^{k} d_{ij} \\
   &e_{ij}=q^{k} e_{ij},\quad e_{ij}=q^{l} e_{ij} \\
 &c_{i-1,j+1}=q^l c_{ij},\quad c_{i+1,j-1}=q^k c_{ij},\\
 &f_{i-1,j+1}=q^k f_{ij},\quad f_{i+1,j-1}=q^l f_{ij},
\end{aligned}
\end{equation*}
 for all $i, j\in\Z_n$.

 If $k+l\equiv 0\ (\mathrm{mod}\ n)$, we obtain that
$$d_{ij}=e_{ij}=0,\quad c_{ij}=q^{ik}c_{0,i+j},\quad f_{ij}=q^{-ik}f_{0,i+j} \quad   \text{for\ all}\ i, j\in\Z_n.  $$
Hence
\begin{equation}\label{eqn2-9}
 \kappa^L(r)=\sum\limits_{i,j=0}^{n-1}q^{ik}c_{0,i+j} u_1x^iy^jz+\sum\limits_{i,j=0}^{n-1}q^{-ik}f_{0,i+j} u_2x^iy^jz.
\end{equation}

Let us now determine the relations between $c_{0,i+j}$ and $f_{0,i+j}$.

By Lemma \ref{lem2-4} and  (\ref{eqn2-9}), we have
 \begin{align*}
&z\ad(\sum\limits_{i,j=0}^{n-1}q^{ik}c_{0,i+j} u_1x^iy^jz)\\
 =&\left(\frac{1}{n}\sum\limits_{r,m=0}^{n-1}q^{-rm}x^r z\cdot u_1 \right)\left(\frac{1}{n}\sum\limits_{g,h=0}^{n-1}q^{-gh}x^g y^mz\left(\sum\limits_{i,j=0}^{n-1}q^{ik}c_{0,i+j} x^iy^jz\right)S(y^{h+m}z) \right)\\
 =&\left(\frac{1}{n}\sum\limits_{r,m=0}^{n-1}q^{r(l-m)}q^{kl} u_2 \right) y^m\left(\frac{1}{n}\sum\limits_{g,h=0}^{n-1}q^{-gh}x^g z\left(\sum\limits_{i,j=0}^{n-1}q^{ik}c_{0,i+j} x^iy^jz\right)zy^{-h} \right)y^{-m}\\
 =&\left(\frac{1}{n}\sum\limits_{m=0}^{n-1}q^{kl}\sum\limits_{r=0}^{n-1}q^{r(l-m)} u_2 \right) \left(\sum\limits_{s,t=0}^{n-1}c_{0,s+t+k}q^{kt} x^{s-m}y^{t+m}z\right)\\
  =&  q^{kl}\left(\sum\limits_{s,t=0}^{n-1}c_{0,s+t+k}q^{kt} u_2x^{s-l}y^{t+l}z\right)\\
   =& \left(\sum\limits_{i,j=0}^{n-1}q^{jk}c_{0,i+j+k}  u_2x^iy^jz\right).
\end{align*}
Similarly,
\begin{align*}
z\ad(\sum\limits_{i,j=0}^{n-1}q^{-ik}f_{0,i+j} u_2x^iy^jz)=&  \left(\sum\limits_{i,j=0}^{n-1}q^{(k-j)k}f_{0,i+j-k} u_1x^iy^jz\right).
\end{align*}
Therefore, the condition that $z\ad \kappa^L(r)=\kappa^L(z\cdot r) $ implies that
\begin{equation*}
 f_{0,i+j}=\pm q^{(i+j)k} c_{0,i+j+k},
\end{equation*}
for all  $i, j\in\Z_n.$

Finally, we get that
 $$\kappa^L(r)=\sum\limits_{i,j=0}^{n-1}q^{ik}c_{0,i+j} u_1x^iy^jz\pm \sum\limits_{i,j=0}^{n-1}q^{jk}c_{0,i+j+k} u_2x^iy^jz.$$

If  $k+l\not\equiv 0\ (\mathrm{mod}\ n)$,  we
can easily get that $\kappa^L(r)=0$ in a similar way. 

Rewriting  $b_{0,i+j}$ and $c_{0,i+j}$  as $b_{i+j}$  and $c_{i+j}$, respectively, we obtain the results.


The proof is finished.
\end{proof}

For the $\mathbb{N}$-graded Koszul $H_{2n^2}$-module algebras $B_u^\pm=T(V)/(I)$, where $V=\mathbb{K}u_1+\mathbb{K}u_2$ and $$I=\{r:=u_1^2\pm q^{kl}u_2^2\}.$$
According to  (\ref{action}), we have
$$x \cdot r = q^{2k}r, \quad y \cdot r = q^{2k}r, \quad z \cdot r = \pm q^{2k^2}r,$$
where $0 \leq k < l \leq n-1$, and $k \equiv l \ (\mathrm{mod}\ \frac{n}{2})$.

The smash products
$B_{u}^{\pm}\sharp\, H_{2n^2}$ can be viewed as the algebras generated by $x, y, z, u_1, u_2$ with the relations
\begin{eqnarray*}\label{eqn1-1}
   \begin{array}{lll}
    &x^n=1,\quad y^n=1, \\
     &xy=yx,\quad zx=yz,\quad zy=xz, \\
    &z^2=\frac{1}{n}\sum\limits_{i,j=0}^{n-1}q^{-ij}x^iy^j,\quad {u_1}^2=\mp q^{kl}{u_2}^2,\\
&xu_1=q^ku_1x, \quad xu_2=q^lu_2x,\\
&yu_1=q^lu_1y, \quad yu_2=q^ku_2y,\\
&zu_1=q^{kl}u_2y^{l}z, \quad zu_2=u_1y^{k}z.\\
   \end{array}
\end{eqnarray*}
As the proof of Theorem \ref{thm2-5},  we can obtain the following results.
\begin{thm}\label{thm2-6} Keep notations as above. Let $B_u^\pm=T(V)/(I)$, where $V=\mathbb{K}u_1+\mathbb{K}u_2$ and $I=\{r:=u_1^2\pm q^{kl}u_2^2\}.$
 Assume that $n$ is even, $0 \leq k< l \leq n-1$,
 and $k\equiv l\ (\mathrm{mod}\ \textstyle\frac{n}{2}),$  the
 deformation
$$\mathscr{U}_{B_{u}^{\pm},\ \kappa} = \frac{\mathbb{K}\langle u_1,u_2 \rangle  \sharp\, H_{2n^2}}{\left( r-\kappa(r)\,|\, {r \in I}\right)},$$
is a nontrivial  $\mathrm{PBW }$-deformation of $B_{u}^{\pm}\sharp\, H_{2n^2}$
if and only if one of the following holds$:$
\begin{enumerate}
  \item [{(a)}] $k=0, \ l=\textstyle\frac{n}{2}$, and
  $$\kappa(r)=
\sum_{i=0}^{n-1}\delta_{+,\pm} a_{ii}x^iy^i
+ \sum_{0 \leq i<j \leq n-1} a_{ij}(x^iy^j \pm x^jy^i)
+ \sum_{i,j=0}^{n-1}\delta_{+,\pm} b_{i+j} x^iy^jz;$$
  \item [{(b)}] $\textstyle\frac{n}{2} $ is even,  $k=\textstyle\frac{n}{4},\ l=\textstyle\frac{3n}{4}$, and
 \begin{align*}
  \kappa(r)=&\pm  \sum\limits_{i,j=0}^{n-1}(-1)^jp^{-\frac{n^2}{4}}
  b_{i+j-\frac{n}{2}}x^iy^jz
    +\sum\limits_{i,j=0}^{n-1}q^{-ik}c_{i+j} u_1x^iy^jz
  \pm\sum\limits_{i,j=0}^{n-1}q^{-jk}c_{i+j-k} u_2x^iy^jz,
  \end{align*}
\end{enumerate}
 where $a_{ij}, b_i, c_i \in \mathbb{K}$ for all $i, j \in \mathbb{Z}_n$, not all of which are zero.
\end{thm}

\section{The PBW-deformations
associated to Koszul dual of AS regular algebras} \label{sec-4}
According to \cite{MYI}, giving a quadratic algebra, we can get its Koszul dual algebra.
This construction retains favorable properties with respect to Koszul duality.

Let $H$ be a Hopf algebra and $A$ an $H$-module algebra.
In this section, it is given the conditions that the Koszul dual  
$A^!$ is  an $H$-module algebra.  In particular, we consider that
  the graded AS regular $H_{2n^2}$-module algebra $A$ of 
  dimension $2$, and  get the conditions that $A^!$ is an $H_{2n^2}$-module algebra.
Then the nontrivial $(0, 1)$-degree PBW-deformations of $A^!\sharp\, H_{2n^2}$ are described.

For a finite-dimensional $\mathbb{K}$-vector space $V$, let $V^\ast = \mathrm{Hom}_{\mathbb{K}}(V, \mathbb{K})$ be the dual space of $V,$ 
and identify $(V \otimes V)^\ast$ with $V^\ast \otimes V^\ast$.
Consider a quadratic  algebra $A = T(V)/\left( I \right)$, where the quadratic  relations $I \subseteq V \otimes V$.
The  Koszul  dual $A^{!}$ of a quadratic $\mathbb{K}$-algebra $A$ is
the quadratic algebra $A^{!}=T(V^*)/\left(I^{\perp}\right),$
where
$$I^{\perp}=\{f \in V^* \otimes V^* \mid f(r)=0\},$$ for all $r \in I$.
For a Koszul algebra $A,$  the dual $A^!$ is also a 
Koszul algebra.

Consider a quadratic algebra
$$A = \frac{\mathbb{K}\langle u_i\,|\, 1\leq i\leq n \rangle}{ \left( r_a |\, 1\leq a\leq m\right)},$$
where $\mathbb{K}\langle u_i\,|\, 1\leq i\leq n\rangle$ denotes the tensor algebra generated by $\{u_i |1\leq i\leq n\}$, and 
$$r_a=\sum_{1\leq i, j\in\leq n}c_a^{ij}u_iu_j \in I,$$
for all $1\leq a\leq m$.
Let $u^s, 1\leq s\leq n$ be the dual basis of $V^*$, that is
$\langle u^s| u_i\rangle=\delta^s_i,$
for all $1\leq i, s\leq n$. Then  the  Koszul dual algebra of $A$ 
can be presented by 
$$A^!=\frac{\mathbb{K}\langle u^s\,|\, 1\leq s\leq n\rangle}{\left(r^b|\, 1\leq b\leq n^2-m\right)} ,$$
where $r^b=\sum\limits_{1\leq s,t\leq n}c^b_{st}u^su^t \in I^{\perp}$, $1\leq b \leq n^2-m$, subject to $\langle r^b|r_a\rangle=0$.

We have the following result.
\begin{prop}\label{lem4-1}
Let $H$ be a Hopf algebra and $A$ a quadratic $H$-module  algebra.
Then $A^{!}$ is a quadratic $H$-module algebra if and only if
\begin{equation}\label{eqn41-1}
 \sum_{(h)} \left(h_{(2)}\cdot u^s\right)\left(h_{(1)}\cdot u^t\right)=\sum_{(h)}\left(h_{(1)}\cdot u^s\right)
 \left(h_{(2)}\cdot u^t\right)
\end{equation}
hold for all $1\leq s, t\leq n$ and $h\in H.$
\end{prop}
\begin{proof}
It is obvious that $A^!$ is also a quadratic algebra. 

Furthermore,  $V^* = \sum_{s=1}^n \mathbb{K} u^s$ is an $H$-module  in the natural 
way:
\begin{equation}\label{h}
(h \cdot u^\ast)(u) = \sum_{(h)} u^\ast(S(h) \cdot u),
\end{equation}
for any $h \in H $, $ u^\ast \in V^*$, and $u \in V$, where $S$ is the antipode of $H$.

We observe that a quadratic algebra $A^!$ is generated by the 
basis $ \{ u^s \mid 1\leq s \leq n \} $ of $V^*$ 
with some quadratic relations. The action on $V^*$ can be naturally extended to the tensor algebra $T(V^*).$  And $T(V^\ast)$ is 
an $H$-module algebra  is equivalent to that
\[
h \cdot (u^s u^t) = \sum_{(h)} (h_{(1)} \cdot u^s) (h_{(2)} \cdot u^t),
\]
for all $1\leq s, t \leq n$.

Now, for all $u^s, u^t \in V^*$ and  $u_i, u_j \in V$, we have
\begin{equation}\label{aln4-2}
\begin{aligned}
 \left(h\cdot \left(u^s u^t\right)\right)\left(u_i u_j\right)
= &\left(u^s u^t\right) \left( S(h)\cdot \left(u_i  u_j\right)\right) \\
=&(u^s u^t) \left(\sum\limits_{(h)}\left(S(h_{(2)})\cdot u_i\right) \left(S(h_{(1)})\cdot u_j\right) \right)\\
=&\sum\limits_{(h)}\left(u^{s}\left(S(h_{(2)}\cdot u_i\right)\right )\left(u^{t}\left(S(h_{(1)}\cdot u_j\right)\right)\\
=&\sum\limits_{(h)}\left(h_{(2)}\cdot u^{s}\right)(u_i) \left(h_{(1)}\cdot u^{t}\right)(u_j)\\
=&\left(\sum\limits_{(h)}\left(h_{(2)}\cdot u^{s}\right) \left(h_{(1)}\cdot u^{t}\right)\right)\left(u_i u_j\right).
\end{aligned}
\end{equation}
Hence,
 $$h\cdot \left(u^su^t\right)=\sum\limits_{h} \left(h_{(1)}\cdot u^s\right)  \left( h_{(2)}\cdot u^t\right)$$ holds if and only if
$$
 \sum_{(h)} \left(h_{(2)}\cdot u^s\right)\left(h_{(1)}\cdot u^t\right)=\sum_{(h)}\left(h_{(1)}\cdot u^s\right)
 \left(h_{(2)}\cdot u^t\right).
$$
Now, we still need to demonstrate that the $H$-actions on $I^\bot$ are closed.
Suppose that
$$r^{b}=\sum\limits_{1\leq s,t\leq n}c^b_{st}u^s u^t\in I^{\perp}.$$
For any relation $r_a=\sum\limits_{1\leq i,j\leq n}c_a^{ij}u_i u_j$ in $I$ and all $h\in H$,
we get that
\begin{align*}
\left(h\cdot r^b\right)(r_a)=&\sum\limits_{1\leq s,t\leq n}c^{b}_{st}\left(\sum\limits_{(h)}\left(h_{(2)}\cdot u^{s}\right) \left(h_{(1)}\cdot u^{t}\right)\right)\left(\sum\limits_{1\leq i,j\leq n}c_a^{ij}u_i u_j\right)\\
= &\left(\sum_{1\leq s,t\leq n}c^{b}_{st}u^s  u^t\right) \left( S(h)\cdot \left(\sum\limits_{1\leq i,j\leq n}c_a^{ij}u_i  u_j\right)\right) \\
=&0.
\end{align*}
Therefore, $h\cdot r^b\in I^\bot$, and  $I^\bot$ is closed 
under the action of $H$.

The result follows.
\end{proof}
 
 Based on the above results, we can get the conditions
 that $\left(A_u^\pm\right)^!$  (resp. $\left(B_u^\pm\right)^!$ ) are 
 quadratic $H_{2n^2}$-module algebras.
 
Recall that  
$${A_u^{\pm}}= \frac{\mathbb{K}\langle u_1, u_2\rangle}{ \left(u_1u_2\pm p^{l^2-k^2}u_2u_1\right)} $$ 
are Koszul $H_{2n^2}$-module algebras.
Let $u^1$ and $u^2$ be the dual basis of $A_1^*$, then 
\begin{alignat*}{3}
    &x\cdot u^1=q^{-k} u^1,
    \quad\quad&&y\cdot u^1=q^{-l}u^1,
    \quad\quad&&z\cdot u^1
    =u^2,\\
    &x\cdot u^2=q^{-l} u^2,
    \quad\quad&&y\cdot u^2=q^{-k} u^2,
    \quad\quad&&z\cdot u^2
    =q^{kl}u_1,
    \end{alignat*}
by  (\ref{action}), and  (\ref{h}).
It is evident that
\[\left(A_u^{\pm}\right)^!= \frac{\mathbb{K}\langle u^1, u^2\rangle}{\left( r^m\,|\,1\leq m\leq 3 \right) },\]
by \cite{MYI},
where $r^1:= u^1u^1=(u^1)^2$, $r^2:= u^2u^2=(u^2)^2$, and $r^3:=u^1u^2\mp p^{k^2-l^2}u^2u^1.$

By   (\ref{eqn41-1}), it is easy to see that
$\left(A_u^{\pm}\right) ^!$ are Koszul $H_{2n^2}$-module algebras
if and only if
$k^2\equiv l^2\ (\mathrm{mod}\ n).$ 
Furthermore, the actions of the elements $x$, $y$, and $z$ on $r^1$, $r^2$ and $r^3$ 
are given by
\begin{equation}\label{a-43}
\begin{array}{lll}
  x\cdot r^1 =q^{-2k}r^1,\quad & y\cdot r^1 =q^{-2l}r^1,  \quad & z\cdot r^1 =q^{kl}r^2,\\
  x\cdot r^2 =q^{-2l}r^2, & y\cdot r^2 =q^{-2k}r^2, & z\cdot r^2 =q^{3kl}r^1,\\
  x\cdot r^3 =q^{-k-l}r^3, & y\cdot r^3 =q^{-k-l}r^3, & z\cdot r^3 = \mp p^{(k+l)^2}r^3.
\end{array}\end{equation}


\begin{thm}\label{thm4-10} Keep notations as above. For  $0\leq k<l \leq n-1$, the deformation
$$\mathscr{U}_{\left(A_u^{\pm}\right)^{!}, \kappa}=\frac{\mathbb{K}\langle u^1,u^2\rangle \sharp\, H_{2n^2}}{\left( r^m-\kappa\left(r^m\right)\,|\, {1\leq m\leq 3} \right)}$$
 is  a nontrivial  $\mathrm{PBW}$-deformation of $(A_u^{\pm})^{!} \sharp\, H_{2n^2}$ if and only if
one of the following holds$:$
\begin{enumerate}
\item[(a)]
$k+l \equiv 0\ (\mathrm{mod}\ n)$, and
\begin{align*}
\kappa\left(r^1\right)&=\sum\limits_{i,j=0}^{n-1}q^{2ik}b_{i+j}x^iy^jz\\
 &  + \sum\limits_{i,j=0}^{n-1}q^{ik}{d_{i+j}} u^1x^iy^jz
  +   \sum\limits_{i,j=0}^{n-1} \delta_{q^{3k},1}{e_{ij}} u^2x^iy^j+\sum\limits_{i,j=0}^{n-1}q^{3ik}{f_{i+j}} u^2x^iy^jz,\\
 \kappa\left(r^2\right)&=  \sum\limits_{i,j=0}^{n-1}q^{(2j+k)k}b_{i+j+2k}x^iy^jz\\
 &+   \sum\limits_{i,j=0}^{n-1} \delta_{q^{3k},1}{e_{ji}} u^1x^iy^j+\sum\limits_{i,j=0}^{n-1}q^{3(j+k)k}{f_{i+j+3k}} u^1x^iy^jz
 +\sum\limits_{i,j=0}^{n-1}q^{jk}{d_{i+j}} u^2x^iy^jz,\\
 \kappa\left(r^3\right)&=\sum_{i=0}^{n-1}{\mathbf\delta}_{+,\pm}\, a'_{ii}x^iy^i
+ \sum_{0 \leq i<j \leq n-1} a'_{ij}(x^iy^j \pm x^jy^i)
+ \sum_{i,j=0}^{n-1}\delta_{+,\pm} c_{i+j} x^iy^jz\\
&+\sum\limits_{i,j=0}^{n-1}q^{ik}b'_{i+j} u^1 x^iy^jz
    \mathbf{\pm}  \sum\limits_{i,j=0}^{n-1}q^{jk}{b'_{i+j+k}} u^2x^iy^jz;
\end{align*}
\item[(b)] $4\mid n$, $k+l\not\equiv 0\ (\mathrm{mod}\ n), \ k+l\equiv 0 \ (\mathrm{mod}\ \textstyle\frac{n}{2})$, and
\begin{align*}
\kappa\left(r^1\right)&=
\sum\limits_{i,j=0}^{n-1}\delta_{q^{k},1}a_{ij}x^iy^j+\sum\limits_{i,j=0}^{n-1}q^{2ik}b_{i+j}x^iy^jz,\\
\kappa\left(r^2\right)&=\sum\limits_{i,j=0}^{n-1}\delta_{q^{k},1}a_{ji}x^iy^j+ \sum\limits_{i,j=0}^{n-1}q^{(2j-l)k}b_{i+j+2k}x^iy^jz,\\
\kappa\left(r^3\right)&=\pm\sum_{i,j=0}^{n-1}(-1)^jp^{-\frac{n^2}{4}}c_{i+j-
\frac{n}{2}}x^iy^jz,
\end{align*}
\end{enumerate}
 where $a_{ij}, a'_{ij}, b_i,  b'_i, c_{i}, d_i, e_{ij}, f_i \in \mathbb{K}$ for all $i, j \in \mathbb{Z}_n$, not all of which are zero.
\end{thm}
\begin{proof} Reviewing the fact  that
 $\left(A_u^{\pm}\right) ^!$ are Koszul $H_{2n^2}$-module algebras
if and only if $k^2\equiv l^2\ (\mathrm{mod}\ n),$
we get that if the statements (a) and (b) hold, $\left(A_u^{\pm}\right)^!$ have 
to be $H_{2n^2}$-module algebras.
And we have the smash product $(A_u^{\pm})^{!} \sharp\, H_{2n^2}.$ 

On the other hand, under the condition that $k^2\equiv l^2\ (\mathrm{mod}\ n)$, if $(A_u^{\pm})^{!} \sharp\, H_{2n^2}$ have nontrivial PBW-deformations, the cases (a) and (b) have to be 
satisfied.

Indeed, it is sufficient to  
consider Theorem \ref{thm1-4}(1)
since $\mathrm{dim}_{\mathbb{K}}\left(A_u^\pm\right)^!_1 = 2$.
This is equivalent to check that
\begin{equation}\label{eqn-yang}
   h\ad\kappa(r^m)=\kappa(h\cdot r^m)
\end{equation} 
for all $1\leq m\leq 3$ and  $h\in H_{2n^2}$ imply that statements (a) and (b).

By Theorem \ref{thm2-5}, 
if
 (\ref{eqn-yang}) for $m=3$ is satisfied, then 
\begin{align*}
\kappa\left(r^3\right)&=\sum_{i=0}^{n-1}{\mathbf\delta}_{+,\pm}\, a'_{ii}x^iy^i
+ \sum_{0 \leq i<j \leq n-1} a'_{ij}(x^iy^j \pm x^jy^i)
+ \sum_{i,j=0}^{n-1}\delta_{+,\pm} c_{i+j} x^iy^jz\\
&+\sum\limits_{i,j=0}^{n-1}q^{ik}b'_{i+j} u^1 x^iy^jz
    \mathbf{\pm}  \sum\limits_{i,j=0}^{n-1}q^{jk}{b'_{i+j+k}} u^2x^iy^jz.
\end{align*}

Let
\begin{align*}
\kappa(r^1) &= \sum_{i,j=0}^{n-1} \left( a_{ij} x^i y^j + b_{ij} x^i y^j z \right), \\
&+ \sum_{i,j=0}^{n-1} \left( c_{ij} u^1 x^i y^j + d_{ij} u^2 x^i y^j z \right) 
 + \sum_{i,j=0}^{n-1} \left( e_{ij} u^1 x^i y^j + f_{ij} u^2 x^i y^j z \right),\\
\kappa(r^2) &= \sum_{i,j=0}^{n-1} \left( a'_{ij} x^i y^j + b'_{ij} x^i y^j z \right), \\
&+ \sum_{i,j=0}^{n-1} \left( c'_{ij} u^1 x^i y^j + d'_{ij} u^2 x^i y^j z \right) 
 + \sum_{i,j=0}^{n-1} \left( e'_{ij} u^1 x^i y^j + f'_{ij} u^2 x^i y^j z \right),
\end{align*}
where  $a_{ij}, a'_{ij}, b_{ij}, b'_{ij}, c_{ij}, c'_{ij}, d_{ij}, d'_{ij}, e_{ij}, e'_{ij}, f_{ij}, f'_{ij} \in \mathbb{K}.$

 By 
  (\ref{a-43}) and the conditions that $x\ad \kappa(r^m)=\kappa(x\cdot r^m)$ and $y\ad \kappa(r^m)=\kappa(y\cdot r^m)$ for all $m\in\{1,2\}$, we obtain that
\begin{equation}\label{a-44}
\begin{array}{llll}
a_{ij}=q^{-2k}a_{ij}, &b_{i-1,j+1}=q^{-2k}b_{ij}, \\
a_{ij}=q^{-2l}a_{ij}, &  b_{i+1,j-1}=q^{-2l}b_{ij},\\
a'_{ij}=q^{-2l}a'_{ij},&b'_{i-1,j+1}=q^{-2l}b'_{ij},\\
a'_{ij}=q^{-2k}a'_{ij},  &
b'_{i+1,j-1}=q^{-2k}b'_{ij},
\end{array}
\end{equation}
and

\begin{equation}\label{b-44}
\begin{array}{llll}
c_{ij}=q^{-k}c_{ij}, &  d_{i-1,j+1}=q^{-k}d_{ij},\\
c_{ij}=q^{-l}c_{ij},& d_{i+1,j-1}=q^{-l}d_{ij},\\
e_{ij}=q^{-2k+l}e_{ij},  &
f_{i-1,j+1}=q^{-2k+l}f_{ij},\\
e_{ij}=q^{-2l+k}e_{ij}, &
f_{i+1,j-1}=q^{-2l+k}f_{ij},\\
c'_{ij}=q^{-2l+k}c'_{ij}, &  d'_{i-1,j+1}=q^{-2l+k}d'_{ij},\\
c'_{ij}=q^{-2k+l}c'_{ij},& d'_{i+1,j-1}=q^{-2k+l}d'_{ij},\\
e'_{ij}=q^{-l}e'_{ij}, &
f'_{i+1,j-1}=q^{-l}f'_{ij},\\
e'_{ij}=q^{-k}e'_{ij},  &
f'_{i-1,j+1}=q^{-k}f'_{ij},\\
\end{array}
\end{equation}
where   $i,j\in \mathbb{Z}_n.$

Let us divide two cases to discuss.
 
{\bf (a).  $k+l\equiv 0\ (\mathrm{mod}\ n)$.}

In this case,  we have $k^2\equiv l^2\ (\mathrm{mod}\ n)$, and $(A_u^{\pm})^{!}$ are $H_{2n^2}$-module algebras.
The smash products  $(A_u^{\pm})^{!}\sharp\, H_{2n^2}$ 
are obtained.

By (\ref{a-43})-(\ref{b-44}), and  the condition that $z\ad \kappa(r^m)=\kappa(z\cdot r^m)$ for all $m\in\{1,2\}$, we obtain that
\begin{equation*}\label{a-45}
\begin{aligned}
&a_{ij}=0,   &b_{ij}=q^{2ik}b_{0,i+j},\\
&a'_{ij}=0,   &b'_{ij}=q^{(2j+k)k}b_{0,i+j+2k},\\
&c_{ij} = 0,   &d_{i,j} = q^{ik}d_{0,i+j}, \\
&e_{ji} = q^{3k}c'_{ij},     &f_{i,j} = q^{3ik}f_{0,i+j}, \\
&c'_{ji} = q^{-3k^2+3k}e_{ij}, &d'_{i,j} = q^{3(j+k)k}f_{0,i+j+3k}, \\
&e'_{ij} = 0,     &f'_{i,j} = q^{jk}d_{0,i+j+k},\\
\end{aligned}
\end{equation*}
where $i,j\in \mathbb{Z}_n.$

Consequently,   $\kappa(r^m)$, $m\in\{1,2\}$  can be written as
\begin{align*}
\kappa(r^1)
&=\sum\limits_{i,j=0}^{n-1}q^{2ik}b_{0,i+j}x^iy^jz
   + \sum\limits_{i,j=0}^{n-1}q^{ik}{d_{0,i+j}} u^1x^iy^jz\\
 & +   \sum\limits_{i,j=0}^{n-1} \delta_{q^{3k},1}{e_{ij}} u^2x^iy^j+\sum\limits_{i,j=0}^{n-1}q^{3ik}{f_{0,i+j}} u^2x^iy^jz,\\
\kappa(r^2)
 &=\sum\limits_{i,j=0}^{n-1}q^{(2j+k)k}b_{0,i+j+2k}x^iy^jz+\sum\limits_{i,j=0}^{n-1}q^{jk}{d_{0,i+j}} u^2x^iy^jz\\
 &+   \sum\limits_{i,j=0}^{n-1} \delta_{q^{3k},1}{e_{ji}} u^1x^iy^j+\sum\limits_{i,j=0}^{n-1}q^{3(j+k)k}{f_{0,i+j+3k}} u^1x^iy^jz,
\end{align*}
were $b_{0,i+j}, d_{0,i+j}, e_{ij}, f_{0,i+j}\in\mathbb{K}.$

{\bf (b).  $k+l\not \equiv 0\ (\mathrm{mod}\ n)$}.

 In this case, to ensure that  
 $(A_u^{\pm})^{!}$ are $H_{2n^2}$-module algebras,
 or equivalently,  that $k^2\equiv l^2\ (\mathrm{mod}\ n)$,
  it is necessary that
 $\textstyle\frac{n}{2}$
  be even and $k+l\equiv 0\ (\mathrm{mod}\ \textstyle{\frac{n}{2}})$. 
 Subsequently, we 
 get the smash products  
 $(A_u^{\pm})^{!}\sharp\, H_{2n^2}.$ 
 
  By (\ref{a-43})-(\ref{b-44}) and  the condition that $z\ad \kappa(r^m)=\kappa(z\cdot r^m)$ for all $m\in \{1,2\}$, we obtain that
\begin{equation*}\label{a-46}
\begin{aligned}
&a_{ij}=\delta_{q^k,1}a_{ij},  \quad  b_{ij}=q^{2ik}b_{0,i+j},\\
&a'_{ij}=\delta_{q^k,1}a'_{ij}, \quad  b'_{ij}=q^{(2j-l)k}b_{0,i+j+2k},\\
&c_{ij}=c'_{ij}=0,  \quad  d_{ij}=d'_{ij}=0,\\
&e_{ij}=e'_{ij}=0,  \quad  f_{ij}=f'_{ij}=0,
\end{aligned}\end{equation*}
where $i,j\in \mathbb{Z}_n$.

Consequently,  $\kappa(r^m)$, $m\in\{1,2\}$  can be written as
\begin{align*}
\kappa(r^1)
&=\sum\limits_{i,j=0}^{n-1}\delta_{q^{k},1}a_{ij}x^iy^j+\sum\limits_{i,j=0}^{n-1}q^{2ik}b_{0,i+j}x^iy^jz,\\
\kappa(r^2)
 & =
 \sum\limits_{i,j=0}^{n-1}\delta_{q^{k},1}a_{ji}x^iy^j+ \sum\limits_{i,j=0}^{n-1}q^{(2j-l)k}b_{0,i+j+2k}x^iy^jz.
\end{align*}

Summarizing the above two cases and rewriting $b_{0,i+j}$,  $d_{0,i+j}$, and $f_{0,i+j}$  as $b_{i+j}$, $d_{i+j}$, and  $f_{i+j}$, respectively, we get the  proof of the  Theorem.
\end{proof}

For the Koszul $H_{2n^2}$-module algebras  
$${B_{u}^{\pm}}= \frac{\mathbb{K}\langle {u_1}, {u_2}\rangle}
{ \left( {u_1}^2\pm q^{kl}{u_2}^2 \right)},$$ 
where $n$ is even, $0\leq k<l\leq n-1$ and $k\equiv l\ (\mathrm{mod}\ \textstyle\frac{n}{2}),$  
let ${u^1}$, ${u^2}$ be dual basis  of $V^*$,  we deduce that
\begin{alignat*}{3}
    &x\cdot {u^1}=q^{-k}{u^1},
    \quad\quad&&y\cdot {u^1}=q^{-l}{u^1},
    \quad\quad&&z\cdot {u^1}
    ={u^2},\\
    &x\cdot {u^2}=q^{-l} {u^2},
    \quad\quad&&y\cdot {u^2}=q^{-k} {u^2},
    \quad\quad&&z\cdot {u^2}
    =q^{kl}{u^1},
    \end{alignat*}
by  (\ref{action}) and  (\ref{h}).
We  also have
\[\left({B_{u}^{\pm}}\right)^ != \frac{\mathbb{K}\langle {u^1}, {u^2}\rangle}{\left( r^m\,|\,1\leq m\leq 3\right)}, \]
where $r^1:= {u^1}{u^2}$, $r^2:= {u^2}{u^1}$, and $r^3:={u^1}{u^1}\mp q^{-kl}{u^2}{u^2}$. Furthermore,
$\left({B_{u}^{\pm }}\right)^!$ are Koszul  $H_{2n^2}$-module algebras only if $k^2\equiv l^2\ (\mathrm{mod}\ n)$ by  (\ref{eqn41-1}). The actions of the elements $x$, $y$, and $z$ on $r^1$, $r^2$, and $r^3$ are given by
\[\begin{array}{lll}
  x\cdot r^1 =q^{-k-l}r^1,\quad & y\cdot r^1 =q^{-k-l}r^1,  \quad & z\cdot r^1 =q^{kl+k^2}r^2,\\
  x\cdot r^2 =q^{-k-l}r^2, & y\cdot r^2 =q^{-k-l}r^2, & z\cdot r^2 =q^{kl+k^2}r^1,\\
  x\cdot r^3 =q^{-2k}r^3, & y\cdot r^3 =q^{-2k}r^3, & z\cdot r^3 = \mp q^{2k^2} r^3.
\end{array}\]

As the proof of Theorem \ref{thm4-10}, we  can get
 the following result.
\begin{thm}\label{thm4-11} Keep notations as above.
Assume that $n$ is even,  $0 \leq k< l \leq n-1$, and
$k\equiv l\ (\mathrm{mod}\ \textstyle\frac{n}{2})$,  the 
deformation
$$\mathscr{U}_{(B_u^{\pm})^{!}, \kappa}=\frac{\mathbb{K}\langle u^1,u^2\rangle \sharp\, H_{2n^2}}{\left( r^m-\kappa\left(r^m\right)\,|\, {1\leq m\leq 3} \right)}$$
 is  a nontrivial  $\mathrm{PBW}$-deformation of $(B_u^{\pm})^{!} \sharp\, H_{2n^2},$ if and only if  one of the following holds$:$
\begin{enumerate}
\item[(a)] $4\mid n$,
$k=\textstyle\frac{n}{4},\ l=\textstyle\frac{3n}{4} $ and
\begin{align*}
\kappa\left(r^1\right)
&=\sum\limits_{i,j=0}^{n-1}a_{ij}x^iy^j+\sum\limits_{i,j=0}^{n-1}b_{i+j}x^iy^jz\\
 &  + \sum\limits_{i,j=0}^{n-1}q^{-ik}{d_{i+j}}u^1x^iy^jz
  + \sum\limits_{i,j=0}^{n-1}q^{ik}{e_{i+j}} u^2x^iy^jz, \\
\kappa\left(r^2\right)
 &=\sum\limits_{i,j=0}^{n-1}a_{ji}x^iy^j+\sum\limits_{i,j=0}^{n-1}b_{i+j}x^iy^jz\\
 &+\sum\limits_{i,j=0}^{n-1}q^{jk}{e_{i+j+k}}u^1x^iy^jz
 +\sum\limits_{i,j=0}^{n-1}q^{(-j+k)k}{d_{i+j-k}} u^2x^iy^jz,\\
\kappa\left(r^3\right)
&= \pm  \sum\limits_{i,j=0}^{n-1}(-1)^jp^{-\frac{n^2}{4}}
  c_{i+j-\frac{n}{2}}x^iy^jz\\
& +\sum\limits_{i,j=0}^{n-1}q^{-ik}b'_{i+j}u^1x^iy^jz
  \pm\sum\limits_{i,j=0}^{n-1}q^{-jk}{b'_{i+j-k}} u^2x^iy^jz;
\end{align*}
\item[(b)]  $4\mid n$, $\ k=0,\ l=\textstyle\frac{n}{2},$  and
\begin{align*}
\kappa\left(r^1\right)
&=\sum\limits_{i,j=0}^{n-1}(-1)^ib_{i+j}x^iy^jz,\\
\kappa\left(r^2\right)&= \sum\limits_{i,j=0}^{n-1}(-1)^{j}b_{i+j+\frac{n}{2}}x^iy^jz,\\
 \kappa\left(r^3\right),
&=\sum_{i=0}^{n-1}\delta_{+,\pm} a'_{ii}x^iy^i
+ \sum_{0 \leq i<j \leq n-1} a'_{ij}(x^iy^j \pm x^jy^i)
+ \sum_{i,j=0}^{n-1}\delta_{+,\pm} c_{i+j} x^iy^jz,
\end{align*}
\end{enumerate}
where 
$a_{ij}, a'_{ij}, b_i, b'_i, c_{i}, d_i, e_{i} \in \mathbb{K}$ for all $i, j \in \mathbb{Z}_n,$ not all of which are zero.
\end{thm}
\begin{proof}
  The proof is   similar to Theorem $\ref{thm4-10}$.
\end{proof}

\section{The  PBW-deformations associated to braided products}\label{sec-3}
For $H_{2n^2}$-module algebras $A_u^\pm$ and $B_u^\pm$, we can form the braided products $A_u^\pm \otimes^c (A_u^\pm)_c^{\mathrm{op}}$ and $B_u^\pm \otimes^c (B_u^\pm)_c^{\mathrm{op}}$, respectively. These products are also $H_{2n^2}$-module algebras.  In general, to describe the PBW-deformations $\mathscr{U}_{A, \kappa}$ for a module algebra $A$ with more  generators and relations, it is  tedious  to calculate $\kappa$  (see Theorem \ref{thm1-4} (2)-(4)). 
In this section, we only explore  nontrivial $0$-degree (i.e. $\kappa^L=0$) PBW-deformations of the smash products of $H_{2n^2}$ with these algebras for avoiding cumbersome calculations.

According to \cite[Lemma 2.10]{PD}, there exist Drinfeld twisting elements in  $H_{2n^2}.$ We choose and fix one of them (others are similar to discuss) as
\[J=\frac{1}{n}\sum\limits_{i,j=0}^{n-1} q^{-ij}x^{j}\otimes y^{i}.\] 
The corresponding braiding is
\begin{equation}\label{eqn3-1}
 c=J\circ \tau,
\end{equation}
where $\tau$ is the flip map.

Let us begin by considering the $H_{2n^2}$-module algebras $A_u^{\pm}$.
\begin{prop}\label{prop3-1} The  $H_{2n^2}$-module algebra structures on ${(A_u^{\pm}})_c^{\mathrm{op}}$, ${A_{u}^{\pm}}\otimes^{c}{(A_u^{\pm})}_c^{\mathrm{op}}$ are given as follows.
\begin{itemize}
  \item[{(a)}] ${(A_u^{\pm}})_c^{\mathrm{op}}$ are  isomorphic to ${A_u^{\pm}}$ as  $H_{2n^2}$-module algebras. 
  \item[{(b)}] As algebras
  $${A_{u}^{\pm}}\otimes^{c}{(A_u^{\pm}})_c^{\mathrm{op}}=\frac{\mathbb{K}\lr{u_1, u_2, v_1, v_2}}{\left( r_m \,|\,  1\leq m\leq 6  \right)} ,$$
  where $0 \leq k<l \leq n-1$, and
   \[\begin{array}{ll}
  r_1 := u_1u_2 \pm p^{l^2-k^2}u_2u_1, \quad & r_2 := v_1v_2 \pm p^{l^2-k^2}v_2v_1, \\
   r_3 := u_1v_1 - q^{-kl}v_1u_1, & r_4 := u_1v_2 - q^{-k^2}v_2u_1, \\
   r_5 := u_2v_1 - q^{-l^2}v_1u_2, & r_6 := u_2v_2 - q^{-kl}v_2u_2.
\end{array}\]
The actions of  $H_{2n^2}$ on  the algebras ${A_{u}^{\pm}}\otimes^{c}{(A_u^{\pm}})_c^{\mathrm{op}}$
are given by   (\ref{action})
 for $\{u_1,u_2\}$ and $\{v_1, v_2\}$ respectively,
 and on the relations in the natural way.
\end{itemize}
\end{prop}
\begin{proof}
{(a)} Let us determine the braided-opposite algebras ${(A_u^{\pm}})_c^{\mathrm{op}}$, which are generated by  $v_1$ and $v_2$. By Definition \ref{defn1-3} (a) and
     the braiding (\ref{eqn3-1}), we obtain that
 \begin{align*}
   &v_1v_1=m_{A_u^{\pm}}\left(\frac{1}{n}\sum\limits_{i,j=0}^{n-1} q^{-ij}x^{j}\otimes y^{i}\right)(u_1\otimes u_1)=q^{kl}u_1u_1,\\
   &v_1v_2=m_{A_u^{\pm}}\left(\frac{1}{n}\sum\limits_{i,j=0}^{n-1} q^{-ij}x^{j}\otimes y^{i}\right)(u_2\otimes u_1)=q^{l^2}u_2u_1,\\
   &v_2v_1=m_{A_u^{\pm}}\left(\frac{1}{n}\sum\limits_{i,j=0}^{n-1} q^{-ij}x^{j}\otimes y^{i}\right)(u_1\otimes u_2)=q^{k^2}u_1u_2,\\
   &v_2v_2=m_{A_u^{\pm}}\left(\frac{1}{n}\sum\limits_{i,j=0}^{n-1} q^{-ij}x^{j}\otimes y^{i}\right)(u_2\otimes u_2)=q^{kl}u_2u_2.
 \end{align*}
These imply that
$$v_1v_2\pm p^{l^2-k^2}v_2v_1 = q^{l^2}u_2u_1 \pm q^{k^2}p^{l^2-k^2}u_1u_2=0.$$
It is clear that $(A_u^{\pm})_c^{\mathrm{op}}$ are isomorphic to $A_u^{\pm}.$

(b) The statement follows in a similar way as in (a). For example, the
 relations $r_3$ and $r_4$ can be deduced by
\begin{align*}
&v_1u_1=\left(m_{A_u^{\pm}}\otimes m_{{(A_u^{\pm}})_c^{\mathrm{op}}}\right)\left(1\otimes c \otimes 1\right)(1\otimes v_1\otimes u_1 \otimes 1)=q^{kl}u_1v_1,\\
&v_2u_1=\left(m_{A_u^{\pm}}\otimes m_{{(A_u^{\pm}})_c^{\mathrm{op}}}\right)\left(1\otimes c \otimes 1\right)(1\otimes v_2\otimes u_1 \otimes 1)=q^{k^2}u_1v_2.
\end{align*}
The proof of $H_{2n^2}$ acting on  the relations of ${A_{u}^{\pm}}\otimes^{c}{(A_u^{\pm}})_c^{\mathrm{op}}$  are also straightforward. For instance, by  (\ref{action}), we obtain
that
\begin{align*}
 z\cdot r_3=&z\cdot \left(u_1v_1-q^{-kl}v_1u_1\right)\\
          =&\frac{1}{n}\sum\limits_{i,j=0}^{n-1}q^{-ij}\left(x^iz\cdot u_1\right)\left(y^{j}z\cdot v_1\right)-q^{-kl}\sum\limits_{i,j=0}^{n-1}q^{-ij}\left(x^iz\cdot v_1)(y^{j}z\cdot u_1\right)\\
          =&\frac{1}{n}q^{2kl}\sum\limits_{i=0}^{n-1}q^{il}\sum\limits_{j=0}^{n-1}q^{j(k-i)}\left(u_2v_2-q^{-kl}v_2u_2\right)\\
          =&q^{3kl}r_6,
\end{align*}
\begin{align*}
 z\cdot r_4=&z\cdot \left(u_1v_2-q^{-k^2}v_2u_1\right)\\
          =&\frac{1}{n}\sum\limits_{i,j=0}^{n-1}q^{-ij}\left(x^iz\cdot u_1\right)\left(y^{j}z\cdot v_2\right)-q^{-k^2}\sum\limits_{i,j=0}^{n-1}q^{-ij}\left(x^iz\cdot v_2\right)\left(y^{j}z\cdot u_1\right)\\
          =&\frac{1}{n}q^{kl}\sum\limits_{i=0}^{n-1}q^{il}\sum\limits_{j=0}^{n-1}q^{j(l-i)}u_2v_1-q^{-k^2}\frac{1}{n}q^{kl}\sum\limits_{i=0}^{n-1}q^{ik}\sum\limits_{j=0}^{n-1}q^{j(k-i)}\\
          =&q^{l(k+l)}\left(u_2v_1-q^{-l^2}v_1u_2\right)\\
          =&q^{l(k+l)}r_5.
\end{align*}
As a result, we get that
\begin{equation}\label{lll}
\begin{array}{lll}
  x\cdot r_1 =q^{k+l}r_1,\quad & y\cdot r_1 =q^{k+l}r_1,  \quad & z\cdot r_1 =\pm p^{(k+l)^2}r_1,\\
  x\cdot r_2 =q^{k+l}r_2, & y\cdot r_2 =q^{k+l}r_2, & z\cdot r_2 =\pm p^{(k+l)^2}r_2,\\
  x\cdot r_3 =q^{2k}r_3, & y\cdot r_3 =q^{2l}r_3, & z\cdot r_3 = q^{3kl}r_6,\\
  x\cdot r_4 =q^{k+l}r_4, & y\cdot r_4 =q^{k+l}r_4, & z\cdot r_4 = q^{l(k+l)}r_5,\\
  x\cdot r_5 =q^{k+l}r_5, & y\cdot r_5 =q^{k+l}r_5, & z\cdot r_5 = q^{k(k+l)}r_4,\\
  x\cdot r_6 =q^{2l}r_6, & y\cdot r_6 =q^{2k}r_6, & z\cdot r_6 = q^{kl}r_3,
\end{array}
\end{equation}
where $0 \leq k<l \leq n-1$.
\end{proof}
 
\begin{thm}\label{thm3-2} 
Keep  notations as above.
For  $0\leq k<l \leq n-1$, the deformation
$$\mathscr{U}_{A_u^{\pm}\otimes^{c}{(A_u^{\pm}})_c^{\mathrm{op}},\ \kappa^C}=\frac{\mathbb{K}\langle u_1,u_2,v_1,v_2\rangle \sharp\, H_{2n^2}}{\left( r_m-\kappa^C(r_m)\, |\, 1\leq m\leq 6\right)}$$
 is a   nontrivial $\mathrm{PBW}$-deformation of $\left(A_u^{\pm}\otimes^{c}{(A_u^{\pm}})_c^{\mathrm{op}}\right) \sharp\, H_{2n^2}$ if and only if one of the following holds$:$
\begin{enumerate}
\item[{(a)}]   $k+l\equiv 0\ (\mathrm{mod}\ n),$
\begin{align*}
\kappa^C(r_1) 
&=\sum\limits_{0 \leq i\leq n-1}\delta_{+,\pm}a_{ii}x^iy^i
+\sum\limits_{\substack{0 \leq i<j \leq n-1,\\ (i-j)k\equiv 0\ (\mathrm{mod}\ n)}}a_{ij}(x^iy^j\pm x^jy^i),\\
\kappa^C(r_2)&=\sum\limits_{0 \leq i\leq n-1}\delta_{+,\pm}b_{ii}x^iy^i
+\sum\limits_{\substack{0 \leq i<j \leq n-1,\\ (i-j)k\equiv 0\ (\mathrm{mod}\ n)}}b_{ij}(x^iy^j\pm x^jy^i),\\
\kappa^C\left(r _3\right)&=\sum\limits_{i,j=0}^{n-1}q^{-2ik}d_{{i+j}}x^iy^jz, \\
\kappa^C\left(r_ 4\right) &=\sum\limits_{\substack {i,j=0,\\ 2(i-j-k)k\equiv 0\ (\mathrm{mod}\ n) }}^{n-1}\delta_{q^{(i-j-k)k},\mp 1}e_{{ij}}x^iy^j-\sum\limits_{i,j=0}^{n-1}q^{-(i+j)k}d_{{i+j-k}}x^iy^jz,\\
\kappa^C\left(r _5\right) &=\sum\limits_{\substack {i,j=0,\\ 2(i-j-k)k\equiv 0\ (\mathrm{mod}\ n) }}^{n-1}\delta_{q^{(i-j-k)k},\mp 1}e_{{ji}}x^iy^j-\sum\limits_{i,j=0}^{n-1}q^{-(i+j)k}d_{{i+j-k}}x^iy^jz,\\
\kappa^C\left(r_ 6\right) &=\sum\limits_{i,j=0}^{n-1}q^{(-2j+3k)k}d_{i+j-2k}x^iy^jz;
\end{align*}
\item[{(b)}] $4\mid n,\ k=0,\ l=\textstyle{\frac{n}{2}}$, and
\begin{align*}
\kappa^C(r_m)&=0,\ \text{for\ all}\ m\in\{1,2\},\\
\kappa^C\left(r _3\right)&= \sum\limits_{i,j=0 }^{n-1}\delta_{(-1)^ip^{\frac{n^2}{4}}, \mp 1}c_{ij}x^iy^j+\sum\limits_{i,j=0}^{n-1}(-1)^{i+j}d_{i+j-\frac{n}{2}}x^iy^jz,\\
 \kappa^C\left(r _4\right) &=\sum\limits_{i,j=0}^{n-1}(-1)^{j+1}d_{i+j}x^iy^jz,\\
 \kappa^C\left(r _5\right)&=\sum\limits_{i,j=0}^{n-1}(-1)^{j+1}d_{i+j}x^iy^jz,\\
 \kappa^C\left(r _6\right)&=\sum\limits_{i,j=0}^{n-1}\delta_{(-1)^ip^{\frac{n^2}{4}}, \mp 1}c_{ji}x^iy^j+\sum\limits_{i,j=0}^{n-1}(-1)^{i+j}d_{i+j-
 \frac{n}{2}}x^iy^jz;
\end{align*}
\item[{(c)}]  $4\mid n,\ (k+l)\equiv 0\ (\mathrm{mod}\ \textstyle{\frac{n}{2}}),
\  k \not\equiv 0\ (\mathrm{mod}\ \textstyle{\frac{n}{2}}),$ and
\begin{align*}
\kappa^C(r_m)&=0, \ \text{for\ all}\ m\in\{1,2\},\\
\kappa^C\left(r _3\right)&=\sum\limits_{i,j=0}^{n-1}q^{-2ik}(-1)^{i+j}d_{i+j-\frac{n}{2}}x^iy^jz,\\
 \kappa^C\left(r _4\right) &=
\sum\limits_{i,j=0}^{n-1}(-1)^{j+1}q^{-(i+j)k}d_{i+j-k}x^iy^jz,\\
\kappa^C\left(r _5\right) &=
\sum\limits_{i,j=0}^{n-1}(-1)^{j+k+1}q^{-(i+j)k}d_{i+j-k}x^iy^jz,\\
\kappa^C\left(r_ 6\right) &=\sum\limits_{i,j=0}^{n-1}(-1)^{i+j+k}q^{(-2j+3k)k}d_{i+j-2k-\frac{n}{2}}x^iy^jz,
\end{align*}
\end{enumerate}
where $a_{ij}, b_{ij}, c_{ij}, d_i, e_{ij}  \in \mathbb{K}$ for all $i, j \in \mathbb{Z}_n$, not all of which are zero.
\end{thm}
\begin{proof}
To get the $0$-th $\mathrm{PBW}$-deformations of $\left(A_u^{\pm}\otimes^{c}{(A_u^{\pm}})_c^{\mathrm{op}}\right) \sharp\, H_{2n^2}$, it suffices to consider
the conditions (1) and (3) of Theorem \ref{thm1-4}(1).
Here, the condition (3) is equivalent to that $\kappa^C$ enjoys 
$\kappa^C \otimes \mathrm{id}=\mathrm{id} \otimes \kappa^C$
on $(I\otimes W)\cap (W\otm I).$

Recall that $W = V \oplus V'$, where $V=\K u_1+\K u_2$, $V'=\K v_1+\K v_2$,
and
${A_{u}^{\pm}}\otimes^{c}{(A_u^{\pm}})_c^{\mathrm{op}}=
T(W)/(I),$
where $I=\left( r_m \, |\, 1\leq m\leq 6\right).$
 Any element in $(I\otimes W) \cap (W \otimes I)$  can be written as
$$c_1 s_1+c_2s_2+c_3s_3+
c_4s_4,$$
where
$$\begin{array}{ll}
s_1:=r_1 v_1-q^{-l^2}r_3 u_2\mp p^{l^2-k^2-2kl}r_5 u_1 & =\quad q^{-l^2-kl}v_1r_1\pm p^{l^2-k^2} u_2r_3+ u_1r_5 , \\
 s_2:=r_1 v_2-q^{-kl}r_4 u_2\mp p^{l^2-3k^2}r_6 u_1 & =\quad q^{-k^2-kl}v_2r_1\pm p^{l^2-k^2} u_2r_4+ u_1r_6 , \\
s_3:=q^{-k^2-kl}r_2 u_1+r_3 v_2\pm p^{l^2-k^2}r_4v_1 & = \quad u_1r_2\mp p^{l^2-3k^2} v_2r_3-q^{-kl} v_1r_4 , \\
s_4:=q^{-l^2-kl}r_2 u_2+r_5 v_2\pm p^{l^2-k^2}r_6v_1 & =\quad u_2r_2\mp p^{l^2-k^2-2kl} v_2r_5-q^{-l^2} v_1r_6,
\end{array}$$
and $c_i\in\K (i=1, 2, 3, 4).$
Also, $s_1, s_2, s_3, s_4$ are linearly independent.
Therefore,  $\{ s_1, s_2, s_3, s_4\}$  
is a basis of $(I\otimes W) \cap (W \otimes I).$  

To ensure that $\mathscr{U}_{A_u^{\pm}\otimes^{c}{(A_u^{\pm}})_c^{\mathrm{op}}}$
is a PBW-deformation,
the following relations must hold for each $1\leq i\leq 4:$
\begin{equation}\label{eqn3-15}
\begin{array}{ll}
\kappa^C(r_1) v_1-q^{-l^2}\kappa^C(r_3) u_2\mp p^{l^2-k^2-2kl}\kappa^C(r_5) u_1 = q^{-l^2-kl}v_1\kappa^C(r_1)\pm p^{l^2-k^2} u_2\kappa^C(r_3)+ u_1\kappa^C(r_5) , \\
\\
 \kappa^C(r_1) v_2-q^{-kl}\kappa^C(r_4) u_2\mp p^{l^2-3k^2}\kappa^C(r_6) u_1= q^{-k^2-kl}v_2\kappa^C(r_1)\pm p^{l^2-k^2} u_2\kappa^C(r_4)+ u_1\kappa^C(r_6) , \\
 \\
q^{-k^2-kl}\kappa^C(r_2) u_1+\kappa^C(r_3) v_2\pm p^{l^2-k^2}\kappa^C(r_4)v_1= u_1\kappa^C(r_2)\mp p^{l^2-3k^2} v_2\kappa^C(r_3)-q^{-kl} v_1\kappa^C(r_4) , \\
\\
q^{-l^2-kl}\kappa^C(r_2) u_2+\kappa^C(r_5) v_2\pm p^{l^2-k^2}\kappa^C(r_6)v_1 = u_2\kappa^C(r_2)\mp p^{l^2-k^2-2kl} v_2\kappa^C(r_5)-q^{-l^2} v_1\kappa^C(r_6).
\end{array}
\end{equation}

Now, we set 
\begin{align*}
&\kappa^C(r_1) = \sum_{i,j=0}^{n-1} (a_{ij}x^iy^j + a'_{ij}x^iy^jz), 
&\kappa^C(r_2) = \sum_{i,j=0}^{n-1} (b_{ij}x^iy^j + b'_{ij}x^iy^jz), \\
&\kappa^C(r_3) = \sum_{i,j=0}^{n-1} (c_{ij}x^iy^j + d_{ij}x^iy^jz), 
&\kappa^C(r_4) = \sum_{i,j=0}^{n-1} (e_{ij}x^iy^j + f_{ij}x^iy^jz), \\
&\kappa^C(r_5) = \sum_{i,j=0}^{n-1} (e'_{ij}x^iy^j + f'_{ij}x^iy^jz),
&\kappa^C(r_6) = \sum_{i,j=0}^{n-1} (c'_{ij}x^iy^j + d'_{ij}x^iy^jz),
\end{align*}
where  $a_{ij}, a'_{ij}, b_{ij}, b'_{ij}, c_{ij}, c'_{ij}, d_{ij}, d'_{ij}, e_{ij}, e'_{ij}, f_{ij}, f'_{ij} \in \mathbb{K}$.

By  (\ref{lll}) and the conditions that $x\ad \kappa^C(r_m)=\kappa^C(x\cdot r_m)$ and $y\ad \kappa^C(r_m)=\kappa^C(y\cdot r_m)$ for all $1\leq m\leq 6,$ we obtain that
\begin{equation}\label{eqn3-2}
\begin{array}{llll}
c_{ij}=q^{2k}c_{ij}, &d_{i-1,j+1}=q^{2k}d_{ij}, \\
c_{ij}=q^{2l}c_{ij}, &  d_{i+1,j-1}=q^{2l}d_{ij},\\
c'_{ij}=q^{2k}c'_{ij},  &
d'_{i+1,j-1}=q^{2k}d'_{i,j},\\
c'_{ij}=q^{2l}c'_{ij},& d'_{i-1,j+1}=q^{2l}d'_{ij},
\end{array}
\end{equation}
and
\begin{equation}\label{eqn3-10}
\begin{aligned}
&a_{ij}=q^{k+l}a_{ij}, \quad\quad\quad \  b_{ij}=q^{k+l}b_{ij}, \quad e_{ij}=q^{k+l}e_{ij},\quad e'_{ij}=q^{k+l}e'_{ij},\\
&a'_{i-1,j+1}=q^{k+l}a'_{ij},\quad b'_{i-1,j+1}=q^{k+l}b'_{ij}, \quad  f_{i-1,j+1}=q^{k+l}f_{ij},\quad f'_{i-1,j+1}=q^{k+l}f'_{ij},\\
  &a'_{i+1,j-1}=q^{k+l}a'_{ij}, \quad b'_{i+1,j-1}=q^{k+l}b'_{ij},\quad f_{i+1,j-1}=q^{k+l}f_{ij}, \quad f'_{i+1,j-1}=q^{k+l}f'_{ij},
\end{aligned}
\end{equation}
where   $i,j\in \mathbb{Z}_n.$

We divide two step to determine the map $\kappa^C$.
  
{\bf (a). $k+l\equiv 0\ (\mathrm{mod}\ n)$.}

By (\ref{eqn3-2})-(\ref{eqn3-10}), and  the conditions that $z\ad \kappa^C(r_m)=\kappa^C(z\cdot r_m)$ for all $1\leq m\leq 6,$ we obtain that
\begin{equation}\label{eqn3-4}
\begin{aligned}
a_{ij} &= \pm a_{ji}, \quad a'_{ij} = \delta_{+,\pm}a'_{0,i+j}, \\
b_{ij} &= \pm b_{ji}, \quad b'_{ij} = \delta_{+,\pm}b'_{0,i+j}, \\
c_{ij} &= 0, \quad\quad \ d_{ij} = q^{-2ik}d_{0,i+j}, \\
e'_{ij} &=e_{ji}, \quad\ f'_{ij} = f_{0,i+j}, \\
c'_{ij} &= 0, \quad\quad f_{ij} = f_{0,i+j}, \\
d'_{ij} &= q^{-2jk+3k^2}d_{0,i+j-2k},
\end{aligned}
\end{equation}
for all $i,j\in \mathbb{Z}_n.$

By  (\ref{eqn3-15}) and (\ref{eqn3-4}), we obtain that
\begin{equation*}\label{ast}
\begin{aligned}
& a_{ij}=q^{(i-j)k}a_{ij}, \quad\ a'_{0,i+j}=0,\\
&b_{ij}=q^{(i-j)k}b_{ij}, \quad\ b'_{0,i+j}=0,\\
&e_{ij}=\mp q^{(i-j-k)k}e_{ij}, \\ & e_{ji}=\mp q^{(i-j+k)k}e_{ji},\\
 &f_{0,i+j}=-q^{-(i+j)k}d_{0,i+j-k},
\end{aligned}
\end{equation*}
for all $i,j\in \mathbb{Z}_n.$

Consequently,
\begin{align*}
\kappa^C(r_1)=&\sum\limits_{0 \leq i\leq n-1}\delta_{+,\pm}a_{ii}x^iy^i+\sum\limits_{\substack{0 \leq i<j \leq n-1,\\ (i-j)k\equiv 0\ (\mathrm{mod}\ n)}}a_{ij}\left(x^iy^j\pm x^jy^i\right),\\
\kappa^C(r_2)=&\sum\limits_{0 \leq i\leq n-1}\delta_{+,\pm}b_{ii}x^iy^i+\sum\limits_{\substack{0 \leq i<j \leq n-1,\\ (i-j)k\equiv 0\ (\mathrm{mod}\ n)}}b_{ij}\left(x^iy^j\pm x^jy^i\right),\\
\kappa^C(r_3)=&\sum\limits_{i,j=0}^{n-1}q^{-2ik}d_{0,i+j}x^iy^jz,\\
\kappa^C(r_4)=&\sum\limits_{{\substack {i,j=0,\\ 2(i-j-k)k\equiv 0\ (\mathrm{mod}\ n) }}}^{n-1}\delta_{q^{(i-j-k)k},\mp 1}e_{ij}x^iy^j- \sum\limits_{i,j=0}^{n-1}q^{-(i+j)k}d_{0,i+j-k}x^iy^jz,\\
\kappa^C(r_5)=&\sum\limits_{\substack {i,j=0,\\ 2(i-j-k)k\equiv 0\ (\mathrm{mod}\ n) }}^{n-1}\delta_{q^{(i-j-k)k},\mp 1}e_{ji}x^iy^j-\sum\limits_{i,j=0}^{n-1}q^{-(i+j)k}d_{0,i+j-k}x^iy^jz,\\
\kappa^C(r_6)=&\sum\limits_{i,j=0}^{n-1}q^{(-2j+3k)k}d_{0,i+j-2k}x^iy^jz,
\end{align*}
where $a_{ij}, b_{ij}, d_{ij}, e_{ij}\in\mathbb{K}.$

{\bf (b). $k+l\not \equiv 0\ (\mathrm{mod}\ n)$}.

 In this case, $n$
 must be even and $k+l\equiv 0\ (\mathrm{mod}\ \textstyle{\frac{n}{2}})$ by  (\ref{eqn3-2})-(\ref{eqn3-10}).

 If $k \not\equiv 0\ (\mathrm{mod}\ \textstyle{\frac{n}{2}})$,
  (\ref{eqn3-2})-(\ref{eqn3-10}) and  the conditions that $z\ad \kappa^C(r^m)=\kappa^C(z\cdot r^m)$ for all $1\leq m\leq 6,$ 
 imply that
\begin{equation}\label{ea3-6}
\begin{aligned}
&a_{ij}=a'_{ij}=b_{ij}=b'_{ij}=c_{ij}=c'_{ij}=e_{ij}=e'_{ij}=0,\\
 &d_{ij}=q^{-2ik}d_{0,i+j}, \quad 
 f_{ij}=(-1)^if_{0,i+j}, \\
 &f'_{ij}=(-1)^{j+l}f_{0,i+j-\frac{n}{2}},
   d'_{ij}=q^{-2jk-3kl}d_{0,i+j-2k},
\end{aligned}
\end{equation}
for all $i,j\in\mathbb{ Z }_n.$

By   (\ref{eqn3-15}) and (\ref{ea3-6}), we deduce that
\begin{equation*}\label{ex}\begin{aligned}
 &d_{0,i+j}=(-1)^{i+j} d_{0,i+j-\frac{n}{2}},\quad
f_{0,i+j}=(-1)^{i+j+1}q^{-(i+j)k}d_{0,i+j-k},
\end{aligned}
\end{equation*}
for all $i,j\in\mathbb{ Z }_n.$

Consequently,  to get the nontrivial PBW-deformations, we must have
 $4\mid n$ and
\begin{align*}
\kappa^C(r_1)=&\kappa^C(r_2)=0,\\
\kappa^C(r_3)=&\sum\limits_{i,j=0}^{n-1}(-1)^{i+j}q^{-2ik}d_{0,i+j-\frac{n}{2}}x^iy^jz,\\
\kappa^C(r_4)=&\sum\limits_{i,j=0}^{n-1}(-1)^{j+1}q^{-(i+j)k}d_{0,i+j-k}x^iy^jz,\\
\kappa^C(r_5)=&\sum\limits_{i,j=0}^{n-1}(-1)^{j+k+1}q^{-(i+j)k}d_{0,i+j-k}x^iy^jz,\\
\kappa^C(r_6)=&\sum\limits_{i,j=0}^{n-1}(-1)^{i+j+k}q^{(-2j+3k)k}d_{0,i+j-2k-\frac{n}{2}}x^iy^jz,
\end{align*}
otherwise, $\kappa^C(r_m)=0$,  for all $1\leq m\leq 6$, and $ d_{ij}\in\mathbb{K}.$

If $k=0,l=\textstyle\frac{n}{2}$,  (\ref{eqn3-2})-(\ref{eqn3-10}) and  the condition that $z\ad \kappa^C(r_m)=\kappa^C(z\cdot r_m)$ for all $1\leq m\leq 6,$ imply that
\begin{equation}\label{ea3-7}
\begin{aligned}
&a_{ij}=b_{ij}=e_{ij}=e'_{ij}=0,  \\ &c'_{ij}=c_{ji}, \quad
  a'_{ij}=0,\quad b'_{ij}=0,\\
  &d_{ij}=d_{0,i+j},  \quad
  f_{ij}=(-1)^if_{0,i+j},\\
  &f'_{ij}=(-1)^{j+1}f_{0,i+j-\frac{n}{2}},\quad
 d'_{ij}=d_{0,i+j},
\end{aligned}
\end{equation}
for all $i,j\in \mathbb{Z}_n.$

By  (\ref{eqn3-15}) and (\ref{ea3-7}), we get that
\begin{equation*}\label{ex}\begin{array}{ll}
(-1)^ic_{ij}=\mp p^{\frac{3n^2}{4}} c_{ij}, & (-1)^jc_{ji}=\mp p^{\frac{-n^2}{4}}  c_{ji},\\
(-1)^ic_{ij}=\mp p^{\frac{n^2}{4}} c_{ij}, &(-1)^jc_{ji}=\mp p^{\frac{-3n^2}{4}}  c_{ji},\\
d_{0,i+j}=(-1)^{i+j}d_{0,i+j-\frac{n}{2}}, &
  f_{0,i+j}=-(1)^{i+j}d_{0,i+j},
\end{array}
\end{equation*}
for all $i,j\in \mathbb{Z}_n.$
Consequently,  to get the nontrivial PBW-deformations, we must have
  $4\mid n,$  and
\begin{align*}
\kappa^C(r_3)=&\sum\limits_{i,j=0}^{n-1}\delta_{(-1)^ip^{\frac{n^2}{4}}, \mp 1}c_{ij}x^iy^j+\sum\limits_{i,j=0}^{n-1}(-1)^{i+j}d_{0,i+j-\frac{n}{2}}x^iy^jz,\\
\kappa^C(r_4)=&\sum\limits_{i,j=0}^{n-1}(-1)^{j+1}d_{0,i+j}x^iy^jz,\\
\kappa^C(r_5)=&\sum\limits_{i,j=0}^{n-1}(-1)^{j+1}d_{0,i+j}x^iy^jz,\\
\kappa^C(r_6)=&\sum\limits_{i,j=0}^{n-1}\delta_{(-1)^ip^{\frac{n^2}{4}}, \mp 1}c_{ji}x^iy^j+\sum\limits_{i,j=0}^{n-1}(-1)^{i+j}
d_{0,i+j-\frac{n}{2}}x^iy^jz.
\end{align*}

Rewriting  $d_{0,i+j}$ as $d_{i}$, we obtain the results.
\end{proof}

For the algebras ${(B_{u}^{\pm})}_c^{\mathrm{op}}$ and ${B_{u}^{\pm}}\otimes^{c}{(B_{u}^{\pm}})_c^{\mathrm{op}},$ we have the following Proposition.

\begin{prop}\label{prop3-2}
The  $H_{2n^2}$-module algebra structures on ${(B_u^{\pm}})_c^{\mathrm{op}}$, ${B_{u}^{\pm}}\otimes^{c}{(B_u^{\pm})}_c^{\mathrm{op}}$ are given as follows.
\begin{itemize}
  \item[{(a)}] $\left({B_u^{\pm}}\right)_c^{\mathrm{op}}$ are  isomorphic to $B_u^{\pm}$ as $H_{2n^2}$-module algebras. 
  \item[(b)]   As algebras
  $${B_{u}^{\pm}}\otimes^{c}{(B_u^{\pm}})_c^{\mathrm{op}}=\frac{\mathbb{K}\lr{u_1, u_2, v_1, v_2}}{\left( r_m \, | \,  1\leq m \leq 6\right)},$$
  where  $0\leq k< l\leq n-1$, $k\equiv l\ (\mathrm{mod}\ \textstyle\frac{n}{2}) $, and
   \[\begin{array}{ll}
  r_1 := {u_1}^2 \pm q^{kl}{u_2}^2, \quad & r_2 := {v_1}^2\pm q^{kl}{v_2}^2, \\
   r_3 := {u_1}{v_1} - q^{-kl}{v_1}{u_1}, & r_4 := {u_1}{v_2} - q^{-k^2}{v_2}{u_1}, \\
   r_5 := {u_2}{v_1} - q^{-l^2}{v_1}{u_2}, & r_6 := {u_2}{v_2} - q^{-kl}{v_2}{u_2}.
\end{array}\]
The actions of  $H_{2n^2}$ on  the algebras ${B_{u}^{\pm}}\otimes^{c}{(B_u^{\pm}})_c^{\mathrm{op}}$
are given by  (\ref{action}) for $\{u_1,u_2\}$ and $\{v_1, v_2\}$ respectively,
 and on the relations in the natural way.
\end{itemize}
\end{prop}
\begin{proof}
 The proof is   similar to  Proposition \ref{prop3-1}.
\end{proof}
 As the proof of Theorem \ref{thm3-2}, we  can get
 the following result.
\begin{thm}\label{thm3-22}
Keep  notations as above. 
 Assume that $n$ is even, $0 \leq k< l \leq n-1$,
 and $k\equiv l\ (\mathrm{mod}\ \textstyle\frac{n}{2})$,  
 the deformation
$$\mathscr{U}_{B_u^{\pm}\otimes^{c}{(B_u^{\pm}})_c^{\mathrm{op}},\ \kappa^C}=\frac{\mathbb{K}\langle u_1,u_2,v_1,v_2\rangle \sharp\, H_{2n^2}}{\left( r_m-\kappa^C(r_m)\, |\,{1\leq m\leq 6}\right)}$$
 is a   nontrivial $\mathrm{PBW}$-deformation of $\left(B_u^{\pm}\otimes^{c}{(B_u^{\pm}})_c^{\mathrm{op}}\right) \sharp\, H_{2n^2}$ if and only if  one of the following holds$:$
\begin{enumerate}
  \item[{(a)}]  $4 \mid  n$, $k=\textstyle\frac{n}{4},$  $l=\textstyle\frac{3n}{4}$, and
\begin{align*}
\kappa^C(r_ m)&=0,\ \text{for\ all}\ m\in\{1,2\},\\
\kappa^C\left(r _3\right)&=\sum\limits_{i,j=0}^{n-1}(-1)^id_{i+j}x^iy^jz, \\
\kappa^C\left(r_4\right) &=   \sum\limits_{\substack{i,j=0,\\ i-j-k \ \mathrm{is} \ \mathrm{even} }}^{n-1}\delta_{q^{(i-j-k)k},-1}e_{ij}x^iy^j\mp\sum\limits_{i,j=0}^{n-1} (-1)^iq^{(j-i)k}d_{i+j+k}x^iy^jz, \\
\kappa^C\left(r_5\right) &=  \sum\limits_{\substack{i,j=0,\\ i-j-k \ \mathrm{is} \ \mathrm{even}   }}^{n-1}\delta_{q^{(i-j-k)k},- 1}e_{ji}x^iy^j\mp\sum\limits_{i,j=0}^{n-1} (-1)^iq^{(j-i)k}d_{i+j+k}x^iy^jz,\\
\kappa^C\left(r_ 6\right)&=\sum\limits_{i,j=0}^{n-1}(-1)^j d_{i+j-\frac{n}{2}}x^iy^jz;
\end{align*}
\item[{(b)}]  $4 \mid n$,  $k=0,\  l=\textstyle\frac{n}{2}$, and
\begin{align*}
 \kappa^C(r_1)
&=\sum\limits_{\substack{0 \leq i\leq n-1,\\ i \ \mathrm{is}\ \mathrm{even}}}\delta_{+,\pm}a_{ii}x^iy^i
+\sum\limits_{\substack{0 \leq i<j \leq n-1, \\ i,j \ \mathrm{is}\ \mathrm{even}}}a_{ij}(x^iy^j\pm x^jy^i),\\
\kappa^C(r_2)&=\sum\limits_{\substack{0 \leq i\leq n-1,\\ i \ \mathrm{is}\ \mathrm{even}}}\delta_{+,\pm}b_{ii}x^iy^i
+\sum\limits_{\substack{0 \leq i<j \leq n-1, \\ i,j \ \mathrm{is}\ \mathrm{even}}}b_{ij}(x^iy^j\pm x^jy^i),\\
\kappa^C\left(r_3\right)&= \sum\limits_{\substack{i,j=0,\\ i,j \ \mathrm{are}\ \mathrm{odd}}}^{n-1}c_{ij}x^iy^j+\sum\limits_{i,j=0}^{n-1}(-1)^{i+j}d_{i+j-\frac{n}{2}}x^iy^jz,\\
\kappa^C\left(r_4\right)&= \sum\limits_{i,j=0}^{n-1}(-1)^{j}d_{i+j}x^iy^jz,\\
\kappa^C\left(r_5\right)&= \sum\limits_{i,j=0}^{n-1}(-1)^{s}d_{i+j-\frac{n}{2}}x^iy^jz,  \\
\kappa^C\left(r_6\right)&= \sum\limits_{\substack{i,j=0,\\ i,j \ \mathrm{are}\ \mathrm{odd}}}^{n-1}c_{ji}x^iy^j+\sum\limits_{i,j=0}^{n-1}(-1)^{i+j}d_{i+j-\frac{n}{2}}x^iy^jz;
\end{align*}
\item[{(c)}] $4\nmid n$,  $k=0,\  l=\textstyle\frac{n}{2}$, and
\begin{align*}
\kappa^C(r_1)
&=\sum\limits_{\substack{0 \leq i\leq n-1,\\ i \ \mathrm{is}\ \mathrm{even}}}\delta_{+,\pm}a_{ii}x^iy^i
+\sum\limits_{\substack{0 \leq i<j \leq n-1, \\ i,j \ \mathrm{is}\ \mathrm{even}}}a_{ij}(x^iy^j\pm x^jy^i),\\
\kappa^C(r_2)&=\sum\limits_{\substack{0 \leq i\leq n-1,\\ i \ \mathrm{is}\ \mathrm{even}}}\delta_{+,\pm}b_{ii}x^iy^i
+\sum\limits_{\substack{0 \leq i<j \leq n-1, \\ i,j \ \mathrm{is}\ \mathrm{even}}}b_{ij}(x^iy^j\pm x^jy^i),\\
\kappa^C\left(r_3\right)&=\sum\limits_{\substack{i,j=0,\\ i,j \ \mathrm{are}\ \mathrm{odd}}}^{n-1}c_{ij}x^iy^j,\\
\kappa^C(r_ m)&=0,\ \text{for\ all}\ m\in\{4,5\},\\
\kappa^C\left(r_6\right)&=\sum\limits_{\substack{i,j=0,\\ i,j \ \mathrm{are}\ \mathrm{odd}}}^{n-1}c_{ji}x^iy^j,
\end{align*}
\end{enumerate}
where $a_{ij}, b_{ij}, c_{ij}, d_{i}, e_{ij}  \in \mathbb{K}$ for all $i, j \in \mathbb{Z}_n$, not all of which are zero.
\end{thm}
\begin{proof}
  The proof is   similar to  Theorem \ref{thm3-2}.
\end{proof}

\section*{ Data Availability Statement}
Data sharing is not applicable to this article, as no datasets were generated or analyzed during the current study.

\section*{ Conflicts of Interest}
The authors declare no conflicts of interest.

\section*{Funding}
The work is supported by National Natural Science Foundation of China (Grant No. 12471038).

\section*{Acknowledgements}
We are deeply grateful for the  comments and suggestions provided by the editors and the reviewers.


\end{document}